\documentclass[11pt,a4paper]{article}
\usepackage[top=2.1cm,bottom=2.1cm,left=2.2cm,right=2.2cm]{geometry}
\usepackage[T1]{fontenc}
\usepackage{lmodern}
\usepackage[dvipsnames,table]{xcolor}
\usepackage{microtype}
\definecolor{accent}{HTML}{1F3A5F}
\definecolor{softblue}{HTML}{EEF3F9}
\usepackage{amsmath,amssymb,amsthm}
\usepackage{booktabs,array}
\usepackage{graphicx}
\usepackage{enumitem}
\usepackage{tikz}
\usetikzlibrary{arrows.meta,positioning,shapes.geometric}
\usepackage{algorithm}
\usepackage{algpseudocode}
\usepackage{subcaption}
\usepackage{float}
\usepackage[font=small,labelsep=colon]{caption}
\usepackage{pgfplots}
\pgfplotsset{compat=1.16}
\usepackage[colorlinks=true,linkcolor=blue!60!black,citecolor=blue!60!black,urlcolor=blue!60!black,pdftitle={A Parallel APSOR--SPIKE Algorithm for Classical and Soft-Layer Obstacle Problems},pdfauthor={Amboru Yalamanda, Peeyush Singh}]{hyperref}
\usepackage{titlesec}
\titleformat{\section}{\Large\bfseries}{\thesection}{1em}{}
\titleformat{\subsection}{\large\bfseries}{\thesubsection}{1em}{}
\titleformat{\subsubsection}{\normalsize\bfseries}{\thesubsubsection}{1em}{}
\titlespacing*{\section}{0pt}{2.4ex plus 0.4ex}{1.4ex}
\titlespacing*{\subsection}{0pt}{1.8ex plus 0.3ex}{0.9ex}

\newtheorem{lemma}{Lemma}
\newtheorem{proposition}{Proposition}
\newtheorem{theorem}{Theorem}
\theoremstyle{remark}
\newtheorem{remark}{Remark}

\newcommand{\nx}{n_x}
\newcommand{\bs}{\texttt}

\algrenewcommand\algorithmicrequire{\textbf{Require:}}
\algrenewcommand\algorithmicensure{\textbf{Ensure:}}

\date{}
\begin{document}
	\pagestyle{plain}
	\begin{center}
		{\LARGE\bfseries A Parallel APSOR--SPIKE Algorithm for\\[4pt] Classical and Soft-Layer Obstacle Problems\par}
		\vspace{22pt}
		{\large Amboru Yalamanda$^{1}$ and Peeyush Singh$^{2}$}\\[8pt]
		{\small $^{1}$ Department of Mathematics and Statistics, VIT-AP UNIVERSITY, Amaravathi, Andhra Pradesh, India\\
			amboru.24phd7073@vitap.ac.in\\[3pt]
			$^{2}$ Department of Mathematics and Statistics, VIT-AP UNIVERSITY, Amaravathi, Andhra Pradesh, India\\
			peeyush.singh@vitap.ac.in}
	\end{center}
	\vspace{10pt}
	\begin{center}\textbf{\small Abstract}\end{center}
	\vspace{-6pt}
	\begin{quote}\small\setlength{\parindent}{0pt}
		This paper presents a parallel scheme for classical obstacle problems and for obstacle problems with a soft layer. The adaptive projected SOR method (APSOR) \cite{MS21} finds the contact set on one process. The linear system defined by the frozen contact set is then solved exactly in parallel by the SPIKE algorithm \cite{PS06}, in which the diagonal blocks are factorized independently and the coupling is condensed into a small reduced system. The remaining nonlinearity of the soft-layer problem \cite{CRWX21} is resolved by a Picard iteration in which every pass reuses the factorization, and a Karush--Kuhn--Tucker (KKT) test stops the iteration. The scheme needs no tuned relaxation parameter. Convergence is proved: APSOR decreases the energy, the SPIKE solve is exact and independent of the number of processes, the Picard iteration is a contraction with factor below $0.052$ on the unit square, and under strict complementarity the loop stops at the unique KKT point. Four classical problems \cite{TSFO14} and four soft-layer examples \cite{CRWX21} in one and two dimensions are computed and compared with exact and analytic solutions (maximum errors $1.1\times10^{-8}$, $7.2\times10^{-14}$ and $6.0\times10^{-4}$), and the KKT test is passed in all eight runs. CPU timings with 1 to 16 MPI ranks show that the SPIKE decomposition takes most of the time in two dimensions, so that the part that is distributed over the ranks is the dominant part of the work.
	\end{quote}
	\vspace{4pt}
	
	\section{Introduction}\label{sec:intro}
	
	The obstacle problem asks for the equilibrium of an elastic membrane that is clamped on the boundary of a domain and held above a rigid obstacle by a load. It is the model variational inequality of the first kind \cite{KS80,Glo84} and occurs in elastostatics, lubrication, torsion and free-boundary problems. After discretization it becomes a convex quadratic problem or, equivalently, a linear complementarity problem with a banded matrix: tridiagonal in one dimension and block tridiagonal with half bandwidth $n$ for the five-point stencil in two dimensions. The difficulty is that the region where the membrane touches the obstacle, the contact set, is not known in advance.
	
	Obstacle problems have to be solved fast and accurately on fine grids, because the contact set and the free boundary are the quantities of interest and their resolution is limited by the mesh. Two difficulties meet here. A sequential sweep method becomes slow on fine grids, and a direct banded solve is exact but costs $O(n^4)$ operations in two dimensions, whereas the contact set that tells a direct solver which equations to keep is exactly what is unknown. Several approaches have been proposed. Tran, Schaeffer, Feldman and Osher \cite{TSFO14} replace the constraint by an exact $L^1$ penalty and minimize without constraints, using split Bregman iterations with a few conjugate-gradient steps for the linear sub-problem. Cheng, Ran, Wang and Xiao \cite{CRWX21} model an obstacle that is a rigid body covered by a deformable soft layer of thickness $g$. The layer pushes back with a normal compliance pressure $p(\psi-u)$, and $g=0$ recovers the classical problem; the nonlinear compliance term has to be resolved on top of the contact set. Projected successive over-relaxation (PSOR) \cite{Cry71} is the classical solver for the discrete problem. Its adaptive version APSOR \cite{MS21} chooses the relaxation parameter by Wolfe conditions, so it needs no tuning, but it remains a sequential Gauss--Seidel sweep.
	
	These methods spend most of their time either in sequential sweeps or in a banded linear solve. SPIKE \cite{PS06,PS07} is a divide-and-conquer banded solver related to the partitioning ideas of Sameh and Kuck \cite{SK78}: the matrix is cut into diagonal blocks that are factorized independently, and the coupling is condensed into a small reduced system. For the obstacle problem the matrix changes whenever the contact set changes, so SPIKE cannot simply be applied to the original problem. Once the contact set is frozen, however, the problem on the free nodes is linear, symmetric and positive definite, and the matrix is a banded matrix to which SPIKE applies. The same idea was used in our earlier work \cite{SY26}, which proposes parallel direct and indirect iterative schemes for general obstacle and free boundary problems. In one dimension the direct scheme is built on the SPIKE banded solver (with a faster recursive variant and a projected version); in two and three dimensions a directional splitting reduces the problem to one-dimensional obstacle problems, and the energy is minimized by a projected gradient iteration with Armijo step selection. Convergence is shown for smooth and nonsmooth convex energies, the direct scheme has a unique solution when the matrix is an M-matrix, and a penalty formulation and an image-deblurring application are also treated. The present paper differs in the algorithm: it combines the adaptive projected SOR predictor with an active-set linearization and a frozen-matrix SPIKE solve, and it treats the soft-layer problem of \cite{CRWX21} as well as the classical one.
	
	In this paper we join the two ideas. APSOR predicts the contact set, SPIKE corrects the solution on that set in parallel, and a KKT test decides whether to stop or to return to APSOR. We give an MPI formulation of this predictor--corrector for classical ($g=0$) and soft-layer ($g>0$) obstacle problems in one and two dimensions, in which the banded Cholesky factorization, the spikes and the reduced system are built once per active set and reused in every Picard pass (Sections~\ref{sec:method}, \ref{sec:how} and~\ref{sec:spike}). We give a self-contained convergence analysis: APSOR decreases the energy for every relaxation parameter below $2/(1+c_\nu h^2/4)$, the SPIKE solve is exact, so that the iterates do not depend on the number of processes, the Picard iteration for the compliance term is a contraction, and the loop terminates at the unique KKT point under strict complementarity (Sections~\ref{sec:theory} and~\ref{sec:conv}). We test the scheme on four classical problems from \cite{TSFO14} and four soft-layer examples from \cite{CRWX21}, including a radial problem with an analytic solution and a non-smooth obstacle (Sections~\ref{sec:problems} and~\ref{sec:results}), and we report CPU timings in two dimensions on grids up to $256\times256$ and in one dimension for matrix orders up to $8192$, with up to 16 MPI ranks and a stage profile that shows where the time goes (Section~\ref{sec:perf}).
	
	The scheme has several advantages. It needs no tuning, because the relaxation parameter of APSOR is chosen by the Wolfe conditions and the only parameters are tolerances, and it treats the classical problem ($g=0$) and the soft-layer problem ($g>0$) with the same code. The linear system is solved exactly by SPIKE, so the accuracy is that of the discretization, up to $7.2\times10^{-14}$ in our tests, and the result does not depend on the number of processes. The factorization, the spikes and the reduced system are computed once per active set and reused in every Picard pass, only one to six factorizations were needed in the two-dimensional runs, and the work and the storage per rank of the factorization fall with the number of ranks (Section~\ref{sec:spike}). Finally, every run ends with a KKT test, so the computed solution comes with a certificate of optimality, and the convergence of the scheme is proved.
	
	The rest of the paper is organized in the following manner. Section~\ref{sec:apsor} introduces APSOR and the active-set linearization, Section~\ref{sec:model} the models, Section~\ref{sec:method} the method, and Section~\ref{sec:how} how its stages work together. Section~\ref{sec:why} gives the reasons for the approach and its advantages, Sections~\ref{sec:theory} and~\ref{sec:conv} the theory, Section~\ref{sec:spike} SPIKE and its parallel implementation, Section~\ref{sec:perf} the performance, Section~\ref{sec:problems} the test problems, Section~\ref{sec:results} the results, and Section~\ref{sec:concl} the conclusions.
	
	\section{APSOR and the Active-Set Linearization}\label{sec:apsor}
	
	Consider the problem
	\begin{equation}\label{eq:min}
		\min_{u\ge \ell}\ J(u),\qquad J(u)=\tfrac12 u^\top A u+\sum_k P(\psi_k-u_k)-f^\top u,
	\end{equation}
	where $A$ is symmetric positive definite, $\ell=\psi-g$ is the lower bound, $p=P'$ is the compliance of the soft layer and $f$ is the load. For $g=0$ the compliance vanishes and \eqref{eq:min} is the classical obstacle problem $\min_{u\ge\psi}\tfrac12u^\top Au-f^\top u$. Writing $G(u)=\nabla J(u)=Au-p(\psi-u)-f$, the first-order conditions are the complementarity system
	\begin{equation}\label{eq:kkt}
		u\ge\ell,\qquad G(u)\ge0,\qquad G_k(u)\,(u_k-\ell_k)=0\quad\text{for every node }k.
	\end{equation}
	
	\paragraph{Stage 1: adaptive projected SOR.} APSOR sweeps the nodes in lexicographic order. For node $(i,j)$ of the two-dimensional grid it forms the Gauss--Seidel value, relaxes it and projects it onto the constraint,
	\begin{equation}\label{eq:apsor}
		\begin{aligned}
			\hat u_{ij}&=\tfrac14\Big(h^2\big(f_{ij}+p(\psi_{ij}-u_{ij})\big)+u_{i-1,j}+u_{i+1,j}+u_{i,j-1}+u_{i,j+1}\Big),\\
			u_{ij}&\leftarrow\max\{(1-\omega)u_{ij}+\omega\hat u_{ij},\ \ell_{ij}\}.
		\end{aligned}
	\end{equation}
	For $g=0$ this is PSOR \cite{Cry71} with a relaxation parameter that changes from sweep to sweep. The parameter is tied to a step size $s$ by $\omega=2s/(2+s)$ and adapted once per sweep from the energy before and after the sweep, following Algorithm~3 of \cite{MS21}. With $d$ the sweep increment and $V$ the energy, the Armijo test $V_{\text{new}}\le V_{\text{old}}+c_1\nabla V\!\cdot d$ and the curvature test $c_2\nabla V\!\cdot d\le\nabla V_{\text{new}}\!\cdot d$ set the next step as in Table~\ref{tab:wolfe}. The constants are $c_1=0.89$ and $c_2=0.95$, the sweep starts at $s=2$ ($\omega=1$), and the step is reset to $s=2$ if $\omega$ leaves $(0.05,1.99)$.
	
	\begin{table}[!htbp]
		\centering\small
		\caption{Step-size rule of APSOR.}\label{tab:wolfe}
		\begin{tabular}{ll}
			\toprule
			Test outcome & Next step size \\
			\midrule
			Armijo holds and curvature holds & $s\leftarrow1.15\,s$ \\
			Armijo holds, curvature fails & $s\leftarrow1.4\,s$ \\
			Armijo fails & $s\leftarrow0.85\,s$ \\
			\bottomrule
		\end{tabular}
	\end{table}
	
	\paragraph{Stage 2: active-set linearization.} At the end of the sweeps, the nodes resting on the lower bound form the active set $\mathcal A=\{k:u_k\le\ell_k\}$ and the other nodes $\mathcal F$ are free. Freezing $\mathcal A$ makes the problem on $\mathcal F$ linear. The modified matrix $\tilde A$ replaces every active row and column by the identity row and column, and the right-hand side collects the known values of the active neighbours,
	\begin{equation}\label{eq:modified}
		\tilde A_{kl}=\begin{cases}\delta_{kl}, & k\in\mathcal A\ \text{or}\ l\in\mathcal A,\\ A_{kl}, & \text{otherwise,}\end{cases}\qquad
		\tilde b_k=\begin{cases}\ell_k, & k\in\mathcal A,\\ f_k+p(\psi_k-u_k)+h^{-2}\sum_{l\in\mathcal A,\ l\sim k}\ell_l, & k\in\mathcal F.\end{cases}
	\end{equation}
	The matrix $\tilde A$ is symmetric positive definite and has the same band as $A$. This is what makes a banded Cholesky factorization, and therefore SPIKE, applicable. In this sense the predictor APSOR turns the nonlinear constraint into a banded linear system, and the corrector SPIKE solves that system exactly.
	
	\section{Mathematical Modeling}\label{sec:model}
	
	\subsection{Obstacle problem}
	Let $\Omega\subset\mathbb R^d$ be bounded and $\varphi$ a given obstacle. The classical problem seeks
	\begin{equation}\label{eq:classical}
		u=\arg\min_{v\ge\varphi}\ \tfrac12\int_\Omega|\nabla v|^2\,dx,\qquad v=0\ \text{on }\partial\Omega,
	\end{equation}
	equivalently $\int_\Omega\nabla u\cdot\nabla(v-u)\,dx\ge0$ for all $v\ge\varphi$. The one-dimensional obstacles of \cite{TSFO14}, Section 6.1, are given on $[0,1]$ and extended symmetrically from $[0,0.5]$ by
	\begin{equation}\label{eq:obs12}
		\begin{aligned}
			\varphi_1(x)&=\begin{cases}100x^2, & 0\le x\le0.25,\\ 100x(1-x)-12.5, & 0.25\le x\le0.5,\end{cases}\\
			\varphi_2(x)&=\begin{cases}10\sin(2\pi x), & 0\le x\le0.25,\\ 5\cos(\pi(4x-1))+5, & 0.25\le x\le0.5,\end{cases}
		\end{aligned}
	\end{equation}
	with exact solutions on $[0,0.5]$
	\begin{equation}\label{eq:exact12}
		u_1(x)=\begin{cases}(100-50\sqrt2)\,x, & 0\le x\le\tfrac1{2\sqrt2},\\ 100x(1-x)-12.5, & \tfrac1{2\sqrt2}\le x\le0.5,\end{cases}\qquad
		u_2(x)=\begin{cases}10\sin(2\pi x), & 0\le x\le0.25,\\ 10, & 0.25\le x\le0.5.\end{cases}
	\end{equation}
	In two dimensions we use the non-smooth obstacle $\varphi_3$ of \cite{TSFO14}, Eq.~(6.4), given in \eqref{eq:phi3}, and a radial obstacle with an analytic solution, both stated in Section~\ref{sec:problems}.
	
	\subsection{Obstacle problem with a soft layer}
	Following \cite{CRWX21}, the obstacle is a rigid body $\psi$ covered by a soft layer of thickness $g\ge0$. The membrane minimizes the energy over $K=\{v\in H_0^1(\Omega):v\ge\psi-g\}$,
	\begin{equation}\label{eq:energy}
		\begin{aligned}
			E(v)&=\tfrac12\int_\Omega|\nabla v|^2\,dx+\int_\Omega P(\psi-v)\,dx-\int_\Omega fv\,dx,\\
			P(r)&=\int_0^rp(s)\,ds,\qquad p(r)=c_\nu\min\{[r]_+,\,g\},
		\end{aligned}
	\end{equation}
	where $c_\nu$ is the stiffness of the layer (we use $c_\nu=1$) and $[r]_+=\max\{r,0\}$. The solution splits the domain into three regions: \emph{no contact} ($u>\psi$), \emph{soft-layer contact} ($\psi-g<u\le\psi$), where the layer pushes back with $p(\psi-u)>0$, and \emph{rigid contact} ($u=\psi-g$), where the membrane is pressed to the bottom of the layer. For $g=0$ the compliance vanishes and \eqref{eq:energy} reduces to \eqref{eq:classical} with a load.
	
	\subsection{Discretization}
	We use $N$ intervals per direction, mesh size $h=1/N$ on the unit square (and $h=4/N$ on $[-2,2]^2$), $n=N-1$ interior nodes per direction and $\nx=n^2$ unknowns in lexicographic order, $k=(i-1)n+j$. The matrix $A$ is the five-point Laplacian scaled by $h^{-2}$ (diagonal $4h^{-2}$, entries $-h^{-2}$ for grid neighbours). It is symmetric positive definite and banded with half bandwidth $n$. In one dimension it is tridiagonal with half bandwidth $1$.

	\subsection{Notation}
	The symbols used in the paper are collected in Table~\ref{tab:notation}. The letter $p$ is used for the compliance of the soft layer and, in the sections on SPIKE and on the timings, for the number of processes; the meaning is fixed by the section.
	\begin{table}[ht]
		\centering\small
		\caption{Notation.}\label{tab:notation}
		\begin{tabular}{p{3.4cm}p{10.2cm}}
			\toprule
			Symbol & Meaning \\
			\midrule
			$\Omega,\ \psi,\ u,\ f$ & domain, obstacle, solution and load of the obstacle problem (Section~\ref{sec:model}) \\
			$g,\ c_\nu$ & thickness and stiffness of the soft layer; $g=0$ is the classical problem \\
			$p(\cdot),\ P(\cdot)$ & compliance pressure $p(r)=c_\nu\min\{[r]_+,g\}$ and its primitive \\
			$\ell=\psi-g$ & lower bound of the discrete problem \\
			$h,\ N,\ n,\ \nx$ & mesh width; number of intervals per direction; number of interior nodes $n=N-1$ per direction; number of unknowns $\nx=n^2$ in two dimensions \\
			$A,\ b$ & discrete Laplacian (symmetric positive definite) and its half bandwidth ($b=1$ in 1D, $b=n$ in 2D) \\
			$J,\ G$ & discrete energy and its gradient; $G_k\ge0$ on active nodes is the multiplier \\
			$\mathcal A,\ \mathcal F$ & active set (nodes on the lower bound) and free nodes \\
			$\tilde A,\ \tilde b$ & matrix and right-hand side of the linearized problem on the free nodes \\
			$\omega,\ s$ & relaxation parameter of APSOR and its step size \\
			$p,\ B_j,\ V_j,\ W_j$ & number of MPI processes; diagonal blocks of $\tilde A$ and their left and right spikes (Section~\ref{sec:spike}) \\
			$q$ & contraction factor of the Picard iteration \\
			tol & tolerances of the stopping tests (Section~\ref{sec:how}) \\
			\bottomrule
		\end{tabular}
	\end{table}
	
	\section{Methodology: APSOR with SPIKE for Obstacle Problems}\label{sec:method}
	
	Algorithm~\ref{alg:main} states the method. It has four stages, which are explained in detail in Sections~\ref{sec:apsor}, \ref{sec:how} and~\ref{sec:spike}: an APSOR predictor on rank~0, an active-set linearization, a SPIKE corrector on all ranks, and a KKT test.
	
	\begin{algorithm}[!htbp]
		\caption{APSOR--SPIKE for the discrete obstacle problem}\label{alg:main}
		\begin{algorithmic}[1]
			\Require $A$, obstacle $\psi$, layer thickness $g$, load $f$, number of processes $p$, tolerances, sweep limit $200{,}000$
			\Ensure $u$ satisfying the KKT conditions \eqref{eq:kkt}
			\State \textbf{Setup:} partition the $\nx$ unknowns into $p$ consecutive blocks of sizes $m_j$ (Section~\ref{sec:spike}); rank~0 sets $\ell=\psi-g$, $u\gets\max(0,\ell)$, $\varepsilon\gets10^{-3}$
			\While{KKT not satisfied}
			\State \textbf{Stage 1 (rank 0):} APSOR sweeps \eqref{eq:apsor} until the increment $h\|u^{n+1}-u^n\|\le\varepsilon$; \bs{MPI\_Bcast} $u$
			\State \textbf{Stage 2 (all ranks):} $\mathcal A\gets\{k:u_k\le\ell_k+10^{-14}\}$; build $A_j,B_j,C_j$ from \eqref{eq:modified}
			\State \textbf{Stage 3 (all ranks):} banded Cholesky $A_j=L_jL_j^\top$; spike tips $V_j,W_j$; \bs{MPI\_Gather} to rank~0; rank~0 factorizes the reduced system \eqref{eq:reduced}
			\For{$q=1,\dots,40$}\hfill$\triangleright$ Picard passes
			\State all ranks: $b_j\gets$ right-hand side \eqref{eq:modified} from the current $u$; $g_j\gets A_j^{-1}b_j$
			\State \bs{MPI\_Gather} the $2b$ end values; rank~0 solves the reduced system; \bs{MPI\_Bcast} $a_{j+1},b_{j-1}$
			\State all ranks: $x_j\gets A_j^{-1}\big(b_j-[0;B_j]a_{j+1}-[C_j;0]b_{j-1}\big)$; \bs{MPI\_Allgatherv} into $u_{\text{new}}$
			\State \textbf{if} $\|u_{\text{new}}-u\|_\infty<10^{-11}$ \textbf{then break}\hfill$\triangleright$ at once if $g=0$; $u\gets u_{\text{new}}$
			\EndFor
			\State \textbf{Stage 4 (rank 0):} test the KKT conditions (Section~\ref{sec:how}); \textbf{if} they fail, $u\gets\max(u,\ell)$, $\varepsilon\gets\max(\varepsilon/10,10^{-6})$
			\EndWhile
			\State \textbf{return} $u$
		\end{algorithmic}
	\end{algorithm}
	
	The matrix $\tilde A$ changes only when the active set changes, so the factorization of Stage~3 is paid once per outer iteration and reused by all Picard passes. The point of the method is that the sequential stage (APSOR) never has to reach a tight tolerance: it only has to identify the contact set, and the parallel exact solve produces the digits.
	
	\section{Stages of the APSOR--SPIKE algorithm}\label{sec:how}
	
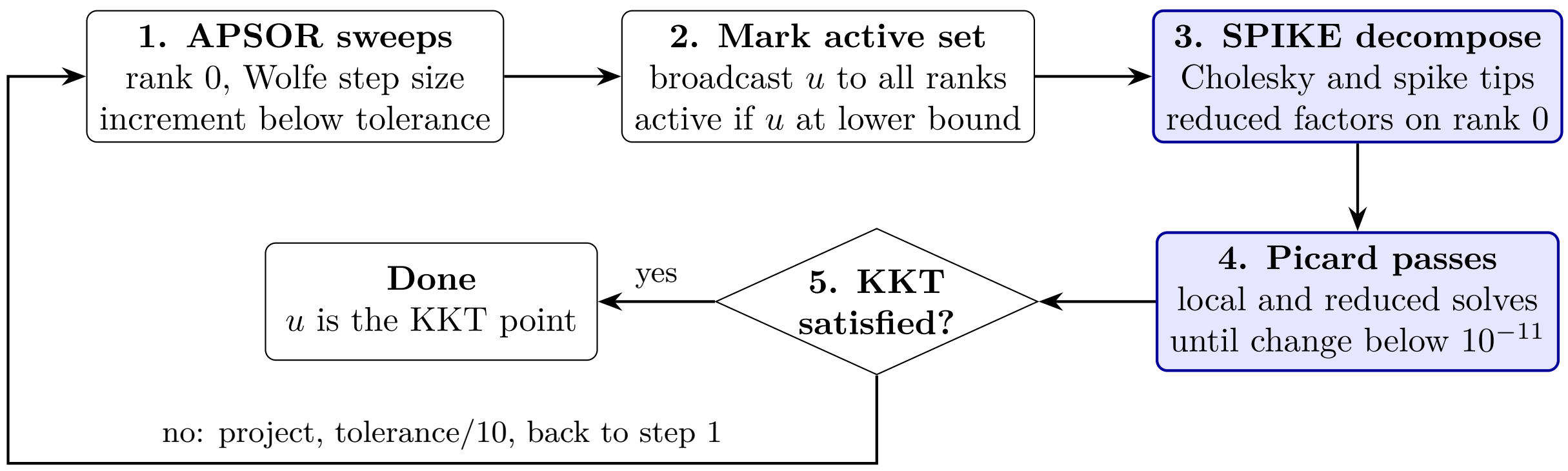
\begin{figure}[!htbp]
		\centering
		\begin{tikzpicture}[
			node distance=9mm and 12mm,
			box/.style={draw, rounded corners=3pt, minimum width=3.4cm, minimum height=1.2cm, align=center, font=\small},
			tint/.style={box, fill=blue!10, draw=blue!60!black, thick},
			dec/.style={draw, diamond, aspect=2.2, align=center, font=\small, inner sep=1pt},
			arr/.style={-{Stealth[length=2.5mm]}, thick}]
			\node[box] (s1) {\textbf{1. APSOR sweeps}\\ rank 0, Wolfe step size\\ increment below tolerance};
			\node[box, right=of s1] (s2) {\textbf{2. Mark active set}\\ broadcast $u$ to all ranks\\ active if $u$ at lower bound};
			\node[tint, right=of s2] (s3) {\textbf{3. SPIKE decompose}\\ Cholesky and spike tips\\ reduced factors on rank 0};
			\node[tint, below=of s3] (s4) {\textbf{4. Picard passes}\\ local and reduced solves\\ until change below $10^{-11}$};
			\node[dec, left=of s4] (s5) {\textbf{5. KKT}\\ \textbf{satisfied?}};
			\node[box, left=of s5] (done) {\textbf{Done}\\ $u$ is the KKT point};
			\draw[arr] (s1)--(s2);
			\draw[arr] (s2)--(s3);
			\draw[arr] (s3)--(s4);
			\draw[arr] (s4)--(s5);
			\draw[arr] (s5)--node[above,font=\footnotesize]{yes}(done);
			\draw[arr] (s5.south) -- ++(0,-0.9) -| node[pos=0.25,above,font=\footnotesize]{no: project, tolerance$/10$, back to step 1} ([xshift=-0.8cm]s1.west |- s1.west) -- (s1.west);
		\end{tikzpicture}
		\caption{Flow of the APSOR--SPIKE algorithm. Shaded boxes are the parallel SPIKE corrector.}
		\label{fig:flow}
	\end{figure}

	Figure~\ref{fig:flow} shows the flow of one run. The subsections below follow its boxes: the predictor (box 1), the active-set marking and corrector (boxes 2 to 4), the Picard iteration (box 4) and the KKT loopback (box 5).
	
	\subsection{Predictor: sequential Gauss--Seidel (APSOR)}
	Rank~0 holds the grid arrays and runs APSOR sweeps from the current iterate, starting from $\max(0,\ell)$ in the first pass. The sweeps stop when the increment is below the stage tolerance, $10^{-3}$ for the first pass. Rank~0 then flattens the grid into the vector $u$ and broadcasts it, and every process marks the active nodes by $u_k\le\ell_k+10^{-14}$. In the pure APSOR mode, which is used as a fallback when memory is short, the tolerance is $10^{-6}$ and the run ends here.
	
	\subsection{Corrector: parallel linearized SPIKE solve}
	\begin{enumerate}[nosep]
		\item \textbf{Constraint locking.} Active rows and columns are replaced by identity rows and columns as in \eqref{eq:modified}. Each process builds its banded block $A_j$ and its coupling blocks $B_j$ and $C_j$ from the stencil without storing a global matrix.
		\item \textbf{Parallel Cholesky.} Every process factorizes $A_j$ independently by banded Cholesky.
		\item \textbf{Spike footprints.} Each process computes the spikes $V_j$ and $W_j$ of Section~\ref{sec:spike} and keeps only their top and bottom $b$ rows. The footprints are gathered on rank~0.
		\item \textbf{Interface solve.} Rank~0 factorizes the reduced system once per active set. In every Picard pass it solves for the interface unknowns from the $2b$ end values of each local solution and broadcasts them.
		\item \textbf{Internal correction.} Each process removes the neighbours' interface contributions from its right-hand side and recovers its block with one banded solve. The blocks are reassembled on every process by \bs{MPI\_Allgatherv}.
	\end{enumerate}
	
	\subsection{Picard iteration for the compliance term}
	When $g>0$ the compliance $p(\psi-u)$ depends on the unknown. With the active set frozen, the code lags it: pass $q$ solves a linear system whose right-hand side uses $u$ from pass $q-1$. The matrix and its SPIKE factorization do not change between passes, so each pass costs one local solve, one reduced solve and one gather. The loop runs at most 40 passes and stops when the largest change in any node is below $10^{-11}$. For $g=0$ the term $p$ vanishes and the first pass is already exact.
	
	\paragraph{Right-hand side and stencil near the active set.} For an active node $b_k=\ell_k$. For a free node $b_k=f_k+p(\psi_k-u_k)$ plus $h^{-2}$ times the lower bound of each active grid neighbour. The left and right neighbours count only if they lie in the same grid row, since at the end of a row the node $k+1$ belongs to the next row and is not a stencil neighbour. The same row test is used when the matrix entries are generated, so the matrix and the right-hand side stay consistent.
	
	\subsection{KKT loopback}
	After the Picard loop, rank~0 evaluates the discrete KKT conditions \eqref{eq:kkt}. The test passes when
	\begin{itemize}[nosep]
		\item $u_k\ge\ell_k-10^{-10}$ at every node (feasibility);
		\item $G_k\ge-10^{-8}$ on active nodes ($u_k\le\ell_k+10^{-12}$);
		\item $\max|G_k|\le10^{-6}$ on free nodes.
	\end{itemize}
	If it passes the run ends. Otherwise $u$ is projected onto $u\ge\ell$, the APSOR tolerance is divided by 10 (never below $10^{-6}$), and the algorithm returns to the predictor. The loop also ends when the total number of APSOR sweeps reaches $200{,}000$. The decision is made on rank~0 and broadcast.

	\paragraph{Execution of the stages.} The two methods are not run at the same time: they are coupled inside the outer loop. Table~\ref{tab:sep} lists, for each stage, where it runs, what it does mathematically and which routine of the code implements it.
	\begin{table}[ht]
		\centering\small
		\caption{Separation of the stages of the algorithm.}\label{tab:sep}
		\begin{tabular}{p{3.2cm}p{2.9cm}p{4.6cm}p{3.2cm}}
			\toprule
			Stage & Execution & Mathematical role & Routine \\
			\midrule
			APSOR sweeps & sequential, rank~0 only & predictor: damps local errors and estimates the active set & \bs{APSOR\_SWEEPS} \\
			Local factorization and spikes & parallel, all ranks & corrector: banded Cholesky of $A_j$, spike tips $V_j,W_j$ & \bs{CHOL\_BAND}, \bs{DECOMPOSE} \\
			Reduced system & rank~0 & interface unknowns between the blocks & \bs{REDUCED\_SOLN} \\
			Picard passes & parallel, all ranks & local solves, retrieval of the blocks, \bs{MPI\_Allgatherv} & \bs{SOLVE\_BAND}, \bs{REST\_PART\_SOLVE} \\
			KKT test & rank~0 & stop or return to APSOR & \bs{KKT\_STATS} \\
			\bottomrule
		\end{tabular}
	\end{table}
	
	\section{Rationale and properties of the approach}\label{sec:why}
	
	This section motivates the combination of APSOR with a SPIKE-based direct solver, compares it with the alternatives for the problems of Section~\ref{sec:model} and states its advantages.
	
	\subsection{APSOR as the active-set predictor}
	The inequality $u\ge\ell$ makes the discrete problem nonlinear and non-smooth. A projected sweep handles it with no matrix at all and no penalty parameter, and the Wolfe-condition rule removes the relaxation parameter that PSOR \cite{Cry71} has to tune for each problem. The energy decreases at every sweep for any $\omega$ below the bound of Proposition~\ref{prop:apsor}. APSOR is quick at identifying which nodes touch the obstacle, because that is a local, high-frequency feature. It is slow at removing the smooth, long-range part of the error, and that is why we do not use it to reach a tight tolerance.
	
	\subsection{SPIKE as the parallel banded solver}
	For a frozen active set the problem on the free nodes is a banded, symmetric positive definite linear system. A serial banded Cholesky factorization of the two-dimensional problem, with $\nx=n^2$ unknowns and half bandwidth $n$, costs $O(n^4)$ operations \cite{GVL13} and is the dominant cost. SPIKE distributes this work over $p$ processes: the diagonal blocks are factorized independently, the coupling is condensed into a reduced system of dimension $2b(p-1)$, and, because the matrix is constant while the active set is, the factorization, the spikes and the reduced system are computed only once per active set. In one dimension the matrix is tridiagonal, a serial Thomas solve costs $O(n)$ and SPIKE cannot beat it at the sizes considered. The one-dimensional problems 6.1 and 6.2 therefore serve to verify the algorithm, and the parallel advantage of SPIKE is in two dimensions, where the half bandwidth is large.
	
	\subsection{Exact versus inexact inner solvers}
	SPIKE with the full spikes and the full reduced system is an \emph{exact} solver. The iterates therefore do not depend on the number of processes or on the partition (Proposition~\ref{prop:spike}), and the convergence theory of Section~\ref{sec:conv} carries over to the parallel code unchanged. This is not true for iterative or domain-decomposition solvers that solve the linear system only approximately, because their inexactness enters the outer iteration and has to be controlled separately.
	
	\subsection{Sequential predictor and parallel corrector}
	A pointwise relaxation sweep is sequential, and information moves through the grid by at most one cell per sweep. A change of the obstacle near the first block therefore needs of the order of $n$ sweeps to reach the last block, so the number of sweeps grows with the grid (Table~\ref{tab:p61} shows $207$ sweeps at $N=512$ and $10\,231$ at $N=8192$ for Problem~6.1). The sweeps are, however, good at what is needed first: they smooth the local error and show which nodes rest on the obstacle. SPIKE complements this. It couples all blocks through the reduced interface system in one solve and does not need a sequence of sweeps to carry the information across the blocks. APSOR is therefore used as a predictor of the active set, and SPIKE as a corrector that gives the exact solution on that set.
	
	The parallel stage has two further properties. The factorization of the diagonal blocks needs no communication, and the data that are exchanged are the spike tips and the interface values, whose size depends on the number of blocks and on the half bandwidth and not on the volume of the global matrix (Section~\ref{sec:spike}). The code is matrix-free and keeps only a banded block per rank, so that the storage per rank is limited. Table~\ref{tab:op} summarizes the two stages.
	\begin{table}[ht]
		\centering\small
		\caption{Comparison of the APSOR predictor and the SPIKE corrector.}\label{tab:op}
		\begin{tabular}{p{3.0cm}p{5.2cm}p{5.2cm}}
			\toprule
			& APSOR predictor & SPIKE corrector \\
			\midrule
			Execution & sequential on rank~0, the other ranks wait & all ranks work on their blocks \\
			Error that it removes & local, high-frequency error; finds the contact set & long-range coupling between the blocks, through the reduced system \\
			Communication & one broadcast of $u$ per run of sweeps & spike tips once per active set; interface values and \bs{MPI\_Allgatherv} per Picard pass \\
			Type & iterative, nonlinear smoother & direct, exact linear solve \\
			Cost in our runs & $99\%$ of the time in 1D at $N=8192$, about $1\%$ in 2D & $84\%$ to $97\%$ of the time in 2D \\
			\bottomrule
		\end{tabular}
	\end{table}
	
	\subsection{Comparison with alternative methods}
	Table~\ref{tab:cmp} summarizes the methods mentioned above. The entries on parallel structure are qualitative.
	
	\begin{table}[!htbp]
		\centering\small
		\caption{Qualitative comparison of solution methods for the discrete obstacle problem.}\label{tab:cmp}
		\begin{tabular}{p{3.3cm}p{3.3cm}p{4.2cm}p{3.6cm}}
			\toprule
			Method & Tuning & Parallel structure & Proven convergence \\
			\midrule
			PSOR \cite{Cry71} & $\omega\in(0,2)$ & sequential sweeps & yes (SPD matrices) \\
			APSOR \cite{MS21} & none (Wolfe rule) & sequential sweeps & Proposition~\ref{prop:apsor} (this paper) \\
			$L^1$-penalty split Bregman \cite{TSFO14} & penalty parameter, shrinkage & CG inner solves, no parallel code & yes, for the penalized problem \\
			APSOR + SPIKE (this paper) & none & factorization once per active set, exact parallel solve & Theorem~\ref{thm:outer}, under strict complementarity \\
			\bottomrule
		\end{tabular}
	\end{table}

	
	\subsection{Strengths of the combined scheme}
	The reasons above give the following advantages of APSOR--SPIKE. The predictor needs no parameter tuning and is cheap, so that the contact set is found with sequential sweeps that cost $O(N)$ each. The linear system for the frozen active set is solved exactly, so that the final accuracy is that of the discretization and not of an iteration. The dominant cost, the banded factorization, is shared among the ranks, and it is paid once per active set and not once per Picard pass. The result is independent of the number of processes (Proposition~\ref{prop:spike}), so that the scheme can be run on any number of ranks without changing the answer. The KKT test certifies the computed solution, and the same code solves the classical and the soft-layer problem.
	
	\section{Mathematical Theory}\label{sec:theory}
	
	\subsection{Well-posedness}
	
	\begin{proposition}[Unique solution]\label{prop:wp}
		Let $A$ be symmetric positive definite and $p$ nondecreasing. Then the discrete problem \eqref{eq:min} has a unique solution $u^*$, and $u^*$ is the unique point satisfying the KKT conditions \eqref{eq:kkt}.
	\end{proposition}
	\begin{proof}
		$P$ is convex because $P'=p$ is nondecreasing, and $u\mapsto P(\psi-u)$ is convex as a convex function of an affine map. With $A$ positive definite, $J$ is strictly convex and coercive. The feasible set $\{u\ge\ell\}$ is closed, convex and nonempty, so a unique minimizer exists. Because the constraints are linear, the KKT conditions are necessary and, by convexity, sufficient.
	\end{proof}
	
	The proposition covers the classical case $g=0$ ($p=0$, a strictly convex quadratic problem) and the soft-layer case $g>0$ with the compliance $p$ of \eqref{eq:energy}.
	
	\subsection{Structure of the matrices}
	
	\begin{lemma}\label{lem:struct}
		(a) Let $\mathcal A$ be any index set. The matrix $\tilde A$ of \eqref{eq:modified} is symmetric positive definite, and so is every diagonal block $A_j$ of any partition of $\tilde A$ into consecutive blocks. (b) On the unit square, the smallest eigenvalue of the five-point matrix is $\lambda_{\min}(A)=\frac{8}{h^2}\sin^2\frac{\pi h}2$, which exceeds $19.4$ for $N\ge8$.
	\end{lemma}
	\begin{proof}
		(a) After a permutation, $\tilde A=\operatorname{diag}(I,A_{\mathcal{FF}})$, where $A_{\mathcal{FF}}$ is a principal submatrix of the positive definite matrix $A$ and hence positive definite. A principal submatrix of a positive definite matrix is positive definite, which gives the claim for the blocks $A_j$. (b) The eigenvalues of the five-point matrix are $\frac4{h^2}\big(\sin^2\frac{i\pi h}2+\sin^2\frac{j\pi h}2\big)$, and the smallest is at $i=j=1$. For $h=1/8$ the value is $19.49$, and it increases as $h$ decreases.
	\end{proof}
	
	\section{Convergence Proof}\label{sec:conv}
	
	We prove convergence of the stages of Algorithm~\ref{alg:main} one at a time and then of the loop. Standing assumptions: $p$ is nondecreasing and Lipschitz with constant $c_\nu$ ($c_\nu=1$ here, and $p=0$ when $g=0$), and $h\le1/8$ on the unit square. Write $J(u)=\tfrac12u^\top Au+\sum_kP(\psi_k-u_k)-f^\top u$, so that $G=\nabla J$.
	
	\subsection{APSOR}
	
	\begin{proposition}[APSOR converges]\label{prop:apsor}
		Suppose every relaxation parameter used by the sweeps lies in $[\omega_{\min},\omega_{\max}]$ with $0<\omega_{\min}$ and $\omega_{\max}<2/(1+c_\nu h^2/4)$. Then $J$ decreases at every sweep, the increments satisfy $\sum_n\|u^{n+1}-u^n\|^2<\infty$, and $u^n\to u^*$.
	\end{proposition}
	\begin{proof}
		Take one node $k$ with current value $u_k$ and gradient $G_k$. Since $h^2G_k/4=u_k-\hat u_k$, the update \eqref{eq:apsor} is $t_k=\max\{u_k-\omega G_k/A_{kk},\ell_k\}-u_k$, a projected coordinate-gradient step of length $\omega/A_{kk}$. The function $\varphi(t)=J(u+te_k)$ has a Lipschitz derivative with constant $L_k=A_{kk}+c_\nu$, because the derivative of $p$ lies in $[0,c_\nu]$. The descent lemma gives
		\[
		J(u+t_ke_k)-J(u)\le G_kt_k+\tfrac{L_k}{2}t_k^2 .
		\]
		Because $u_k$ is feasible and $u_k+t_k$ is the projection of $u_k-\omega G_k/A_{kk}$ onto $[\ell_k,\infty)$, the projection inequality gives $G_kt_k\le-(A_{kk}/\omega)t_k^2$. Combining the two,
		\[
		J(u+t_ke_k)-J(u)\le -\frac{A_{kk}}{2}\Big(\frac{2}{\omega}-1-\frac{c_\nu}{A_{kk}}\Big)t_k^2\le -c_*\,t_k^2,\qquad c_*>0 .
		\]
		The constant $c_*$ is positive exactly when $\omega<2/(1+c_\nu h^2/4)$. On the unit square with $h\le1/8$ this bound is at least $1.992$, so the safeguard $\omega<1.99$ in the code is enough. Each node changes once per sweep, so summing over the sweep gives $J(u^n)-J(u^{n+1})\ge c_*\|u^{n+1}-u^n\|^2$. $J$ is bounded below, hence the increments are square-summable and tend to $0$, and the iterates stay in a bounded set because $J$ is coercive. Let $\bar u$ be a limit point. The intermediate iterates inside a sweep differ from $u^n$ by at most the sweep increment, which tends to $0$, so each coordinate step vanishes in the limit. By continuity of $G$ and of the projection, $\bar u_k=\max\{\bar u_k-\bar\omega G_k(\bar u)/A_{kk},\ell_k\}$ for every $k$, which means $G_k(\bar u)\ge0$ where $\bar u_k=\ell_k$ and $G_k(\bar u)=0$ where $\bar u_k>\ell_k$. This is the KKT system, so $\bar u=u^*$ by Proposition~\ref{prop:wp}. A bounded sequence with a single limit point converges.
	\end{proof}
	
	The Wolfe-condition rule only chooses $\omega$ inside the safeguard interval, so the proof does not depend on it. The increment used as the stopping test tends to zero by the proposition.
	
	\subsection{SPIKE and independence of the number of processes}
	
	\begin{proposition}[The SPIKE solve is exact]\label{prop:spike}
		Let $\tilde A$ be the modified matrix of \eqref{eq:modified}, partitioned into $p$ consecutive diagonal blocks $A_j$ (of any sizes at least $b$, the half bandwidth) with coupling blocks $B_j$ and $C_j$. Then the reduced interface system \eqref{eq:reduced} is nonsingular, the SPIKE solution equals $\tilde A^{-1}\tilde b$ for every $\tilde b$, and the block recursion of Section~\ref{sec:spike} is well defined. Consequently the Picard iterates, and the whole outer iteration, do not depend on $p$ or on the partition, in exact arithmetic.
	\end{proposition}
	\begin{proof}
		By Lemma~\ref{lem:struct}, every $A_j$ is positive definite, so $D=\operatorname{diag}(A_1,\dots,A_p)$ and the spikes exist, and $S=D^{-1}\tilde A$ is nonsingular because $\tilde A$ is. The system $Sx=g$ with $g=D^{-1}\tilde b$ consists of the reduced equations for the top and bottom $b$ rows of each block and of the retrieval equations for all rows. If the reduced system had a nontrivial null vector $y$, extending it by the retrieval equations with $g=0$ would give a nontrivial null vector of $S$, a contradiction. So the reduced system has a unique solution and the retrieval equations return the unique solution of $Sx=g$, that is of $\tilde A x=\tilde b$. The recursion eliminates $b_0,\dots,b_{j-1}$ in turn and leaves $P_ja_j=r_j$. The same argument applied to the leading blocks $0,\dots,j$, a principal submatrix and thus again positive definite, shows that this system has a unique solution for every right-hand side, hence $P_j$ is nonsingular. The map $u\mapsto u_{\text{new}}$ of a Picard pass depends on the linear solve only through its exact result, and the APSOR stage runs on rank~0 only. So no iterate depends on $p$.
	\end{proof}
	
	\subsection{Picard iteration}
	
	\begin{proposition}[Picard iteration is a contraction]\label{prop:picard}
		For a fixed active set, let $T(u_{\mathcal F})=A_{\mathcal{FF}}^{-1}\big(f_{\mathcal F}+p(\psi_{\mathcal F}-u_{\mathcal F})+c\big)$, where $c$ collects the lower bounds of the active neighbours. If $q=c_\nu/\lambda_{\min}(A)<1$, then $\|T(u)-T(v)\|\le q\|u-v\|$ and the Picard iteration converges linearly to the unique fixed point.
	\end{proposition}
	\begin{proof}
		By Cauchy interlacing, $\lambda_{\min}(A_{\mathcal{FF}})\ge\lambda_{\min}(A)$, so $\|A_{\mathcal{FF}}^{-1}\|\le1/\lambda_{\min}(A)$. Since $p$ is Lipschitz with constant $c_\nu$, $\|p(\psi-u)-p(\psi-v)\|\le c_\nu\|u-v\|$. Banach's fixed-point theorem finishes the argument.
	\end{proof}
	
	On the unit square, Lemma~\ref{lem:struct}(b) gives $q<0.052$ for $c_\nu=1$ and $N\ge8$. Because the stopping test uses the largest change and $\|\cdot\|_\infty\le\|\cdot\|_2$, the change after $k$ passes is at most $q^k$ times the first change, so a first change of order $1$ reaches $10^{-11}$ in about nine passes, well inside the cap of 40. The fixed point satisfies $G_k=0$ on all free nodes.
	
	\subsection{The outer loop}
	
	\begin{theorem}[The loop terminates at the KKT point]\label{thm:outer}
		Assume strict complementarity at $u^*$: $G_k(u^*)>0$ where $u^*_k=\ell_k$, and $u^*_k>\ell_k$ elsewhere. Then there is a sweep index $K$ such that, for $n\ge K$, the active set found by APSOR equals the true active set $\mathcal A^*$. Consequently, once the APSOR stage stops at a sweep $n\ge K$, the SPIKE stage returns $u^*$ up to the Picard tolerance and rounding, and the KKT test passes.
	\end{theorem}
	\begin{proof}[Proof (sketch)]
		By Proposition~\ref{prop:apsor}, $u^n\to u^*$. For a node in $\mathcal A^*$, $G_k(u^n)$ stays above some $\delta>0$ and $u^n_k-\ell_k\to0$, so the projected update $u_k-\omega G_k/A_{kk}$ falls below $\ell_k$ and the node is set exactly to $\ell_k$, after which it stays there. For a node outside $\mathcal A^*$, $u^n_k>\ell_k+\delta/2$ for large $n$, so it is never active. Hence the active set equals $\mathcal A^*$ for $n\ge K$. The tolerance is reduced by a factor 10 after every failed test, and the sweep count at which the increment drops below it tends to infinity unless the iterate is already $u^*$, so the stage eventually stops at $n\ge K$, provided that $K$ is reached before the tolerance hits its floor $10^{-6}$. Then $u^*$ satisfies the linearized system of Section~\ref{sec:apsor} for $\mathcal A^*$: it has $u^*_{\mathcal A}=\ell$ and $G_{\mathcal F}(u^*)=0$ on the free nodes. By Propositions~\ref{prop:spike} and~\ref{prop:picard} that system has a unique solution, found by SPIKE, so it equals $u^*$. The KKT tolerances are positive, so the test passes.
	\end{proof}
	
	\begin{remark}[Interpretation]
		(1) The theorem is a statement about the exact outer loop and gives no iteration count; the counts that were measured are in Section~\ref{sec:results}. (2) Without strict complementarity, or if the index $K$ lies beyond the point where the tolerance reaches its floor $10^{-6}$, Proposition~\ref{prop:apsor} still holds and the cap of $200{,}000$ APSOR sweeps ends the loop in every case. (3) No rate is proved for APSOR. (4) The analysis covers the discrete problem only. The discretization error with respect to the continuous solution is not analysed here; Cheng et al.\ \cite{CRWX21} prove an $O(h)$ estimate for a finite element scheme under $H^2$ regularity, which does not cover our finite-difference scheme.
	\end{remark}
	
	\section{The SPIKE Algorithm and its Parallel Implementation}\label{sec:spike}
	
	\subsection{The algorithm}
	We summarize \cite{PS06} for the system $\tilde Ax=\tilde b$ with $\tilde A$ banded. Partition $\tilde A$ into $p$ diagonal blocks $A_j$ of order $m_j$ with coupling blocks $B_j$ (to the following block) and $C_j$ (to the preceding block), each of size $b\times b$. Write $\tilde A=DS$ with $D=\operatorname{diag}(A_1,\dots,A_p)$. The matrix $S$ has identity diagonal blocks and the \emph{spikes}
	\begin{equation}\label{eq:spikes}
		V_j=A_j^{-1}\begin{bmatrix}0\\B_j\end{bmatrix},\qquad W_j=A_j^{-1}\begin{bmatrix}C_j\\0\end{bmatrix}.
	\end{equation}
	First solve $Dg=\tilde b$ in parallel. Then, with $a_j$ and $b_j$ the first and last $b$ entries of the solution block $x_j$, and $Vt_j,Vb_j,Wt_j,Wb_j$ and $gt_j,gb_j$ the top and bottom parts of the spikes and of $g_j$, the interface unknowns satisfy the block tridiagonal reduced system
	\begin{equation}\label{eq:reduced}
		a_j+Vt_j\,a_{j+1}+Wt_j\,b_{j-1}=gt_j,\qquad b_j+Vb_j\,a_{j+1}+Wb_j\,b_{j-1}=gb_j,
	\end{equation}
	with $V_p=W_1=0$, of dimension $2b(p-1)$ after eliminating the unknowns that do not appear. Rank~0 solves it by block elimination, writing $b_j=c_j-F_ja_{j+1}$, $P_j=I-Wt_jF_{j-1}$, $T_j=P_j^{-1}Vt_j$, $S_j=Wb_jF_{j-1}$ and $F_j=Vb_j+S_jT_j$. The matrices $P_j$ are LU-factorized once per active set and every Picard pass only reuses them. Once the interface values are known, every rank recovers its block independently,
	\begin{equation}\label{eq:retrieve}
		x_j=A_j^{-1}\big(b_j-[0;B_j]\,a_{j+1}-[C_j;0]\,b_{j-1}\big).
	\end{equation}
	Our implementation builds the spikes in full, keeps only their top and bottom $b$ rows ($2b^2$ numbers per spike) and solves the full reduced system, so it is an exact solver (Proposition~\ref{prop:spike}). The truncated reduced system, which is block diagonal for diagonally dominant matrices \cite{PS06}, is not used.
	
	\subsection{MPI structure}
	The code is Fortran 90 with MPI \cite{MPI15}. Rank $j$ owns $m_j$ consecutive unknowns. One outer iteration performs the following steps.
	\begin{enumerate}[nosep]
		\item Rank~0 runs APSOR and broadcasts $u$ (\bs{MPI\_Bcast}).
		\item Every rank forms the active set, builds and factorizes $A_j$ and computes the spike tips \eqref{eq:spikes}; the tips are gathered on rank~0, which factorizes the reduced system.
		\item In every Picard pass each rank solves $A_jg_j=b_j$ and sends the $2b$ interface entries of $g_j$ to rank~0 (\bs{MPI\_Gather}).
		\item Rank~0 solves the reduced system and broadcasts the $2b(p-1)$ interface values.
		\item Each rank recovers its block by \eqref{eq:retrieve} and all blocks are assembled on every rank by \bs{MPI\_Allgatherv}.
	\end{enumerate}
	In two dimensions a row-wise partition needs every block to contain at least two grid rows, so $p\le(N-1)/2$. In one dimension the runs have $b=1$.
	
\begin{figure}[!htbp]
		\centering
		\begin{tikzpicture}[blk/.style={draw,minimum height=0.8cm,font=\footnotesize\ttfamily,inner sep=3pt}]
			\node[blk,minimum width=3cm] (r0) at (0,0) {Rank 0 data};
			\node[blk,minimum width=3.4cm,right=0pt of r0] (r1) {Rank 1 data};
			\node[blk,minimum width=1.4cm,right=0pt of r1] (rd) {$\cdots$};
			\node[blk,minimum width=3.6cm,right=0pt of rd] (rp) {Rank $p{-}1$ data};
			\draw[-{Stealth}] ([yshift=-8mm]r0.south west)--(r0.south west);
			\draw[-{Stealth}] ([yshift=-8mm]r1.south west)--(r1.south west);
			\draw[-{Stealth}] ([yshift=-8mm]rp.south west)--(rp.south west);
			\node[font=\footnotesize\ttfamily,anchor=north west] at ([yshift=-8mm]r0.south west) {displs(0)=0};
			\node[font=\footnotesize\ttfamily,anchor=north west] at ([yshift=-8mm,xshift=-3mm]r1.south west) {displs(1)=cnts(0)};
			\node[font=\footnotesize\ttfamily,anchor=north west] at ([yshift=-8mm,xshift=-3mm]rp.south west) {displs($p{-}1$)};
		\end{tikzpicture}
		\caption{Variable-displacement layout of the global vector \bs{unew}. Block $i$ starts at \bs{displs(i)} and has length \bs{cnts(i)}.}
		\label{fig:displs}
	\end{figure}
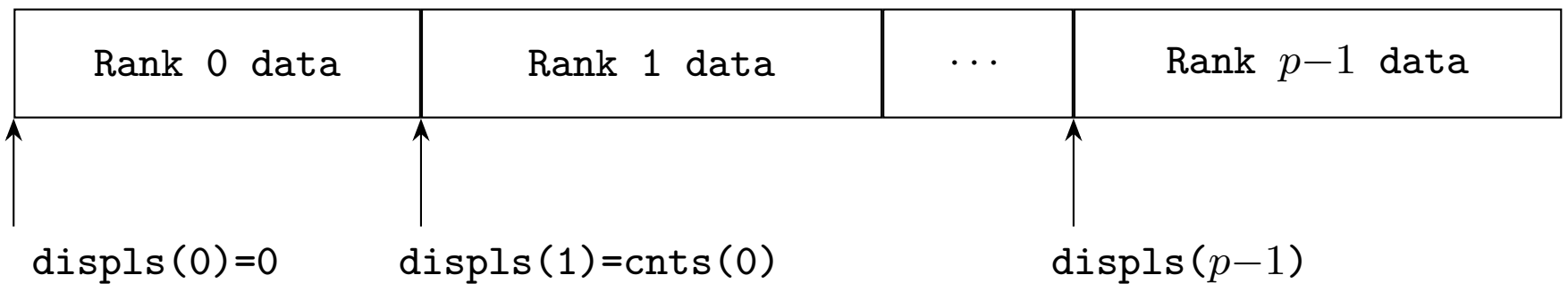

	\paragraph{Variable block sizes.} Since $\nx=(N-1)^2$ is in general not divisible by $p$, block sizes differ by at most one and the uniform \bs{MPI\_Allgather} cannot be used. Two integer arrays of length $p$ describe the layout on every rank (the counts and displacements of \eqref{eq:cnts}),
	\begin{equation}\label{eq:cnts}
		\mathrm{cnts}(i)=\begin{cases}\lfloor \nx/p\rfloor+1, & i<\operatorname{mod}(\nx,p),\\ \lfloor \nx/p\rfloor, & \text{otherwise,}\end{cases}\qquad
		\mathrm{displs}(i)=i\,\lfloor \nx/p\rfloor+\min\{i,\operatorname{mod}(\nx,p)\},
	\end{equation}
	for $i=0,\dots,p-1$. The call is
	\begin{center}\small\bs{MPI\_Allgatherv(xloc, m, MPI\_DOUBLE\_PRECISION,}\\\bs{\ \ unew, cnts, displs, MPI\_DOUBLE\_PRECISION, MPI\_COMM\_WORLD, ierr)}\end{center}
	\noindent Each rank sends its segment \bs{xloc} of length $m$ and receives the whole vector in one collective step, so no staging loop is needed. Figure~\ref{fig:displs} shows the layout. Block $i$ couples only to blocks $i\pm1$, a chain topology; the interface exchange itself goes through rank~0.
	
	\subsection{Cost model and memory}
	Let $m$ be the block size (about $\nx/p$), $b$ the half bandwidth ($n$ in 2D, $1$ in 1D). Table~\ref{tab:cost} lists the leading-order costs and where they are paid.
	
	\begin{table}[!htbp]
		\centering\small
		\caption{Leading-order cost of the stages (model, not measured).}\label{tab:cost}
		\begin{tabular}{lll}
			\toprule
			Operation & Where & Cost \\
			\midrule
			APSOR sweep & rank 0 & $O(\nx)$ per sweep \\
			Banded Cholesky of $A_j$ & every rank & $O(mb^2)$ per active set \\
			Spike tips $V_j,W_j$ & every rank & $O(mb^2)$ per active set \\
			Reduced system factorization & rank 0 & $O(pb^3)$ per active set \\
			Local solve & every rank & $O(mb)$ per Picard pass \\
			Reduced system solve & rank 0 & $O(pb^2)$ per Picard pass \\
			\bottomrule
		\end{tabular}
	\end{table}
	
	In 2D with $m=n^2/p$ and $b=n$, the per-rank work grows as $n^4/p$, while the serial reduced factorization on rank~0 grows as $pn^3$ and increases with $p$. Communication has two parts. Once per active set each rank sends its spike tips ($2b^2$ numbers per spike). In each Picard pass only vectors move: $2b$ numbers gathered per rank, $2b(p-1)$ numbers broadcast, and $m$ numbers per rank in the \bs{MPI\_Allgatherv}, so $O(bp+\nx)$ numbers in total. These estimates are ours and are not measured; they explain the qualitative behaviour of Section~\ref{sec:perf} but are not fitted to it.
	
	The code is matrix-free: a matrix entry is computed on demand from the stencil and the active-set flags. Each rank keeps its block in banded storage that is overwritten by the Cholesky factor, plus the spike tips; rank~0 additionally stores the reduced factors. Before the SPIKE stage the code estimates the memory needed on rank~0 as
	\begin{equation}\label{eq:mem}
		\text{mem (GB)}=\frac{8\left[(b+1)\,m_{\max}+10\,b^2+12\,b^2\,p+12\,\nx+10\,n^2\right]}{10^{9}},
	\end{equation}
	and falls back to pure APSOR if the estimate \eqref{eq:mem} exceeds a limit (default 6~GB). For $N=256$ and $p=16$ the formula gives about $0.12$~GB.
	
	\section{Numerical Performance}\label{sec:perf}
	
	\subsection{Measurements}
	Table~\ref{tab:time} gives the CPU time of rank~0 of the MPI program for the two-dimensional problem (Example 1, $g=0.1$) on grids $N\times N$ with $N=16,\dots,256$ and $p=1,2,4,8,16$ processes, with $p\le(N-1)/2$; a dash means that the run was not made. The one-dimensional timings are in Section~\ref{sec:perf1d}. Rank~0 carries the APSOR stage and the reduced solve, so it is the slowest rank (Table~\ref{tab:profile}). The stage profile of Table~\ref{tab:profile} splits the time into the APSOR stage, the SPIKE decomposition (Cholesky, spikes, reduced factorization) and the Picard passes. Figure~\ref{fig:timing} plots Table~\ref{tab:time} (right panel) and Table~\ref{tab:time1d} (left panel).
	
	\begin{table}[!htbp]
		\centering\small
		\caption{2D obstacle problem: CPU time of rank~0 (s), $N\times N$ grid, $p$ processes. The problem is Example 1 (Section~\ref{sec:ex1grids}).}\label{tab:time}
		\begin{tabular}{rccccc}
			\toprule
			$N$ & $p=1$ & $p=2$ & $p=4$ & $p=8$ & $p=16$ \\
			\midrule
			16 & 0.0041 & 0.0037 & 0.0046 & -- & -- \\
			32 & 0.0185 & 0.0249 & 0.0239 & 0.0243 & -- \\
			64 & 0.178 & 0.296 & 0.221 & 0.274 & 0.340 \\
			128 & 1.558 & 3.036 & 2.131 & 2.286 & 3.148 \\
			256 & 40.94 & 115.23 & 83.64 & 100.66 & 58.60 \\
			\bottomrule
		\end{tabular}
	\end{table}
	
	\paragraph{Speed-up and efficiency.} With $T_p$ the CPU time of rank~0 for $p$ processes, the speed-up is $S_p=T_1/T_p$ and the efficiency is $E_p=S_p/p$. Table~\ref{tab:speed2d} gives both for the two-dimensional runs, and Figure~\ref{fig:speedeff} plots them against the ideal values $S_p=p$ and $E_p=1$. The largest speed-up is $S_2=1.13$ ($E_2=56\%$) on the $16\times16$ grid; for $N\ge32$ the speed-up is between $0.36$ and $0.81$, and the efficiency is $18$--$37\%$ at $p=2$ and $3$--$4\%$ at $p=16$. On the largest grid, $N=256$, the runs give $S_2=0.36$, $S_4=0.49$, $S_8=0.41$ and $S_{16}=0.70$, so that $16$ processes are the fastest multi-process run there.
	
	\begin{table}[!htbp]
		\centering\small
		\caption{2D obstacle problem (Example 1): speed-up $S_p=T_1/T_p$ and efficiency $E_p=S_p/p$ of rank~0, computed from Table~\ref{tab:time}.}\label{tab:speed2d}
		\begin{tabular}{r|cc|cc|cc|cc|cc}
			\toprule
			& \multicolumn{2}{c|}{$p=1$} & \multicolumn{2}{c|}{$p=2$} & \multicolumn{2}{c|}{$p=4$} & \multicolumn{2}{c|}{$p=8$} & \multicolumn{2}{c}{$p=16$} \\
			$N$ & $S_p$ & $E_p$ & $S_p$ & $E_p$ & $S_p$ & $E_p$ & $S_p$ & $E_p$ & $S_p$ & $E_p$ \\
			\midrule
			16 & 1.00 & 100\% & 1.13 & 56\% & 0.90 & 23\% & -- & -- & -- & -- \\
			32 & 1.00 & 100\% & 0.74 & 37\% & 0.77 & 19\% & 0.76 & 10\% & -- & -- \\
			64 & 1.00 & 100\% & 0.60 & 30\% & 0.81 & 20\% & 0.65 & 8\% & 0.52 & 3\% \\
			128 & 1.00 & 100\% & 0.51 & 26\% & 0.73 & 18\% & 0.68 & 9\% & 0.49 & 3\% \\
			256 & 1.00 & 100\% & 0.36 & 18\% & 0.49 & 12\% & 0.41 & 5\% & 0.70 & 4\% \\
			\bottomrule
		\end{tabular}
	\end{table}

\begin{figure}[tbph!]
	\centering
	\includegraphics[width=\linewidth]{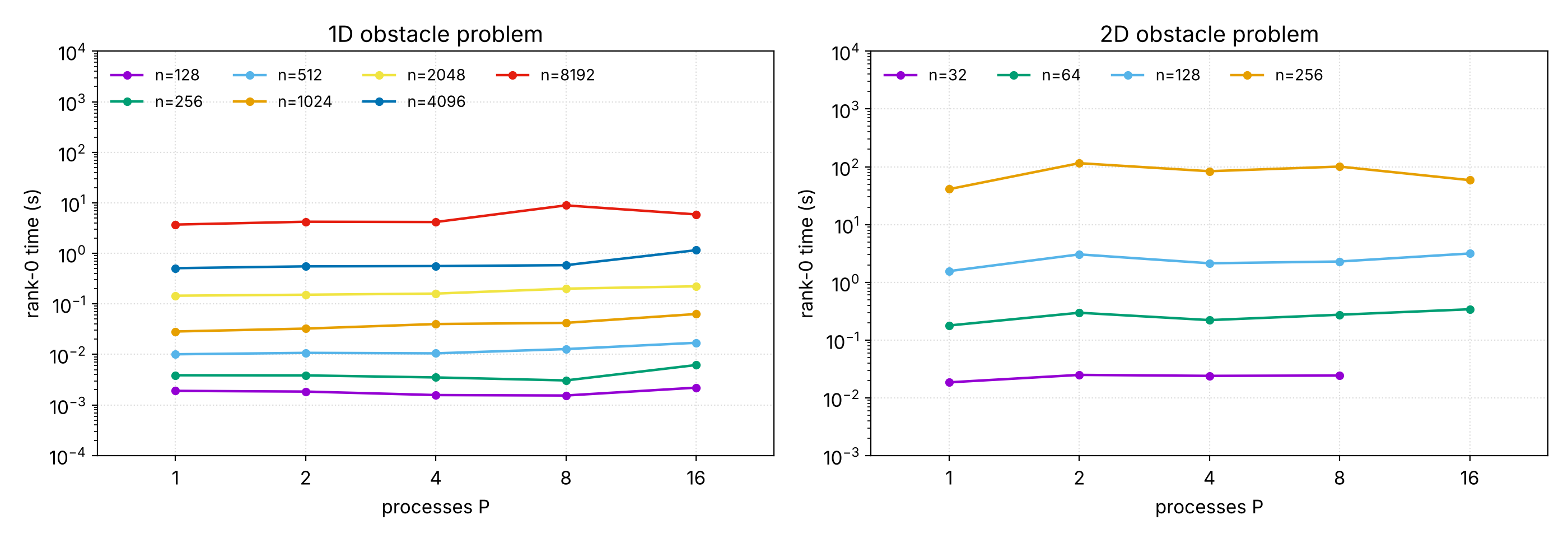}
	\caption{CPU time of rank~0 against the number of MPI processes $p$ (log--log axes). Left: one-dimensional obstacle problem (Problem 6.1), from Table~\ref{tab:time1d}. Right: two-dimensional soft-layer problem (Example 1, $g=0.1$), from Table~\ref{tab:time}.}
	\label{fig:figure03rank0time1d2d}
\end{figure}
	
	\begin{table}[!htbp]
		\centering\small
		\caption{Stage profile (s). ``Workers'' is the CPU time of the non-root ranks. Profile times have a $0.01$~s resolution. The rows are the run outputs of Example 1.}\label{tab:profile}
		\begin{tabular}{rrrrrrrr}
			\toprule
			$N$ & $p$ & Rank 0 & Workers & APSOR & Decompose & Picard & Profile total \\
			\midrule
			64 & 1 & 0.178 & -- & 0.01 & 0.09 & 0.06 & 0.18 \\
			64 & 2 & 0.296 & 0.280 & 0.01 & 0.23 & 0.03 & 0.30 \\
			64 & 4 & 0.221 & 0.205 & 0.01 & 0.18 & 0.02 & 0.22 \\
			64 & 8 & 0.274 & 0.258 & 0.01 & 0.22 & 0.02 & 0.27 \\
			64 & 16 & 0.340 & 0.324 & 0.02 & 0.27 & 0.03 & 0.34 \\
			\midrule
			128 & 1 & 1.558 & -- & 0.04 & 1.14 & 0.36 & 1.56 \\
			128 & 2 & 3.036 & 3.022 & 0.04 & 2.80 & 0.18 & 3.04 \\
			128 & 4 & 2.131 & 2.117 & 0.04 & 1.98 & 0.10 & 2.13 \\
			128 & 8 & 2.286 & 2.271 & 0.07 & 2.09 & 0.11 & 2.29 \\
			128 & 16 & 3.148 & 3.122 & 0.15 & 2.79 & 0.18 & 3.15 \\
			\midrule
			256 & 1 & 40.94 & -- & 0.39 & 34.27 & 6.20 & 40.94 \\
			256 & 2 & 115.23 & 115.16 & 0.44 & 111.22 & 3.48 & 115.23 \\
			256 & 4 & 83.64 & 83.57 & 0.45 & 81.34 & 1.77 & 83.64 \\
			256 & 8 & 100.66 & 100.59 & 0.44 & 98.02 & 2.11 & 100.66 \\
			256 & 16 & 58.60 & 58.53 & 0.97 & 56.12 & 1.40 & 58.60 \\
			\bottomrule
		\end{tabular}
	\end{table}
	
	\subsection{Observations}
	\textbf{Growth with the grid.} With $p=1$, doubling $N$ multiplies the time by $4.5$, $9.7$, $8.7$ and $26.3$ (from $N=16$ to $256$). A banded factorization costs $O(n^4)$, a factor of $16$ per doubling, so the first factors are below this and the last is above it. Part of the growth is the number of factorizations (outer APSOR to SPIKE iterations) of the runs, which was $2$, $2$, $3$, $4$ and $5$ for $N=16,\dots,256$ (Table~\ref{tab:ex1}). Dividing the decomposition time by this number gives $0.03$, $0.29$ and $6.9$~s per factorization at $N=64$, $128$ and $256$, factors of $9.5$ and $24$; the last again exceeds $16$. 
	
	\textbf{Where the time goes.} With one process at $N=256$, the SPIKE decomposition takes $34.27$~s of $40.94$~s ($84\%$), the Picard and reduced solves $6.20$~s and APSOR $0.39$~s ($1\%$). With two processes the decomposition takes $111.22$~s of $115.23$~s, about $97\%$. In every row of Table~\ref{tab:profile} the decomposition is the largest stage and APSOR stays below $0.5$~s. This confirms the cost model: the sequential predictor is cheap and the cost is in the factorization that is paid once per active set.
	
	\textbf{Balance.} The ranks are well balanced. At $N=256$ and $p=16$ the workers finish within $0.01$~s of each other and rank~0 is only $0.07$~s behind them. The extra work on rank~0 (APSOR and the reduced solve) is therefore small.

	\subsection{One-dimensional problems}\label{sec:perf1d}
	Table~\ref{tab:time1d} gives the rank-0 CPU time of the one-dimensional solver (bandwidth $b=1$, tridiagonal matrix of order $N$) for $N=16,\dots,8192$ and $p=1,\dots,16$ processes. The problem timed is Problem~6.1 (Section~\ref{sec:pindep}); each run also reports the CPU time of every rank.
	
	\begin{table}[!htbp]
		\centering\small
		\caption{1D obstacle problem: CPU time of rank~0 (s) for matrix order $N$ and $p$ processes.}\label{tab:time1d}
		\begin{tabular}{rccccc}
			\toprule
			$N$ & $p=1$ & $p=2$ & $p=4$ & $p=8$ & $p=16$ \\
			\midrule
			16 & $6.11\times10^{-4}$ & $7.24\times10^{-4}$ & $5.44\times10^{-4}$ & -- & -- \\
			32 & $8.51\times10^{-4}$ & $9.20\times10^{-4}$ & $9.80\times10^{-4}$ & $8.50\times10^{-4}$ & -- \\
			64 & $1.52\times10^{-3}$ & $1.39\times10^{-3}$ & $1.41\times10^{-3}$ & $1.03\times10^{-3}$ & $2.04\times10^{-3}$ \\
			128 & $1.90\times10^{-3}$ & $1.84\times10^{-3}$ & $1.57\times10^{-3}$ & $1.54\times10^{-3}$ & $2.20\times10^{-3}$ \\
			256 & $3.87\times10^{-3}$ & $3.85\times10^{-3}$ & $3.51\times10^{-3}$ & $3.06\times10^{-3}$ & $6.17\times10^{-3}$ \\
			512 & 0.0101 & 0.0107 & 0.0105 & 0.0128 & 0.0170 \\
			1024 & 0.0284 & 0.0325 & 0.0400 & 0.0421 & 0.0630 \\
			2048 & 0.145 & 0.152 & 0.160 & 0.200 & 0.223 \\
			4096 & 0.509 & 0.553 & 0.560 & 0.584 & 1.155 \\
			8192 & 3.712 & 4.226 & 4.172 & 8.907 & 5.896 \\
			\bottomrule
		\end{tabular}
	\end{table}
	
	\paragraph{Speed-up and efficiency in one dimension.} Table~\ref{tab:speed1d} gives $S_p$ and $E_p$ for Problem~6.1, which Figure~\ref{fig:speedeff} plots together with the two-dimensional values. The largest speed-up is $S_8=1.48$ ($E_8=18\%$) at $N=64$; for $N\ge512$ the speed-up is between $0.42$ and $0.96$ for every $p>1$, and the efficiency at $p=16$ is $3$--$5\%$ for $N\ge128$. For $N\le256$ the times are a few milliseconds, so that the differences from $1$ are of the size of the run-to-run variation of such short runs.
	
	\begin{table}[!htbp]
		\centering\small
		\caption{1D obstacle problem (Problem 6.1): speed-up $S_p=T_1/T_p$ and efficiency $E_p=S_p/p$ of rank~0, computed from Table~\ref{tab:time1d}.}\label{tab:speed1d}
		\begin{tabular}{r|cc|cc|cc|cc|cc}
			\toprule
			& \multicolumn{2}{c|}{$p=1$} & \multicolumn{2}{c|}{$p=2$} & \multicolumn{2}{c|}{$p=4$} & \multicolumn{2}{c|}{$p=8$} & \multicolumn{2}{c}{$p=16$} \\
			$N$ & $S_p$ & $E_p$ & $S_p$ & $E_p$ & $S_p$ & $E_p$ & $S_p$ & $E_p$ & $S_p$ & $E_p$ \\
			\midrule
			16 & 1.00 & 100\% & 0.84 & 42\% & 1.12 & 28\% & -- & -- & -- & -- \\
			32 & 1.00 & 100\% & 0.93 & 46\% & 0.87 & 22\% & 1.00 & 13\% & -- & -- \\
			64 & 1.00 & 100\% & 1.09 & 55\% & 1.08 & 27\% & 1.48 & 18\% & 0.74 & 5\% \\
			128 & 1.00 & 100\% & 1.03 & 52\% & 1.21 & 30\% & 1.24 & 15\% & 0.86 & 5\% \\
			256 & 1.00 & 100\% & 1.01 & 50\% & 1.10 & 28\% & 1.27 & 16\% & 0.63 & 4\% \\
			512 & 1.00 & 100\% & 0.94 & 47\% & 0.96 & 24\% & 0.79 & 10\% & 0.59 & 4\% \\
			1024 & 1.00 & 100\% & 0.88 & 44\% & 0.71 & 18\% & 0.67 & 8\% & 0.45 & 3\% \\
			2048 & 1.00 & 100\% & 0.96 & 48\% & 0.91 & 23\% & 0.73 & 9\% & 0.65 & 4\% \\
			4096 & 1.00 & 100\% & 0.92 & 46\% & 0.91 & 23\% & 0.87 & 11\% & 0.44 & 3\% \\
			8192 & 1.00 & 100\% & 0.88 & 44\% & 0.89 & 22\% & 0.42 & 5\% & 0.63 & 4\% \\
			\bottomrule
		\end{tabular}
	\end{table}
	
\begin{figure}[tbph!]
	\centering
	\includegraphics[width=0.92\linewidth]{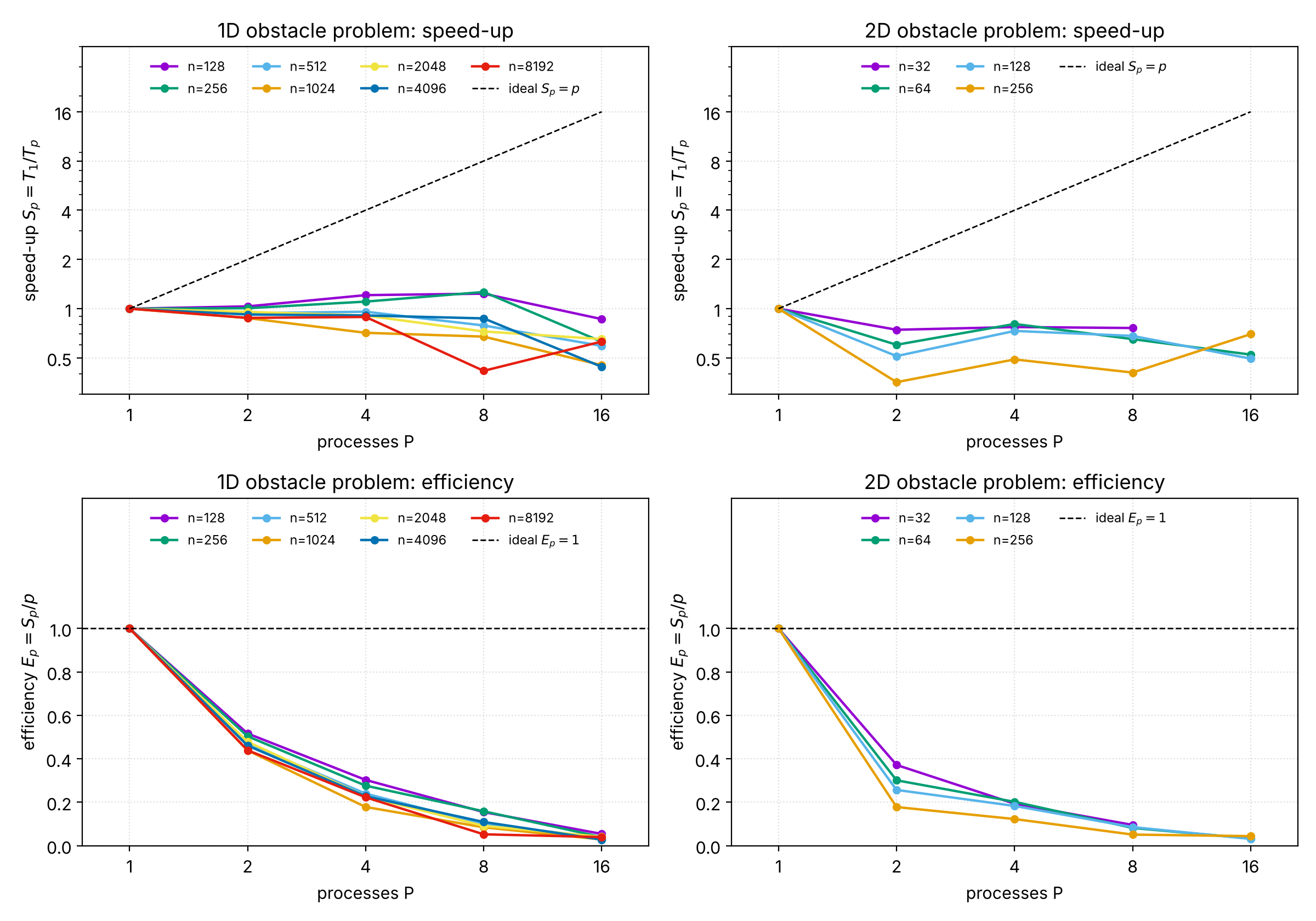}
	\caption{Speed-up $S_p=T_1/T_p$ (top, log scale) and efficiency $E_p=S_p/p$ (bottom) against the number of MPI processes $P$, for the one-dimensional obstacle problem (left, Table~\ref{tab:speed1d}) and the two-dimensional soft-layer problem (right, Table~\ref{tab:speed2d}). The dashed lines are the ideal values $S_p=p$ and $E_p=1$.}
	\label{fig:figure04speedupefficiency1d2d}
\end{figure}

	\paragraph{Observations.} \textbf{Cost structure in one dimension.} A serial tridiagonal solve costs $O(N)$ and SPIKE with $b=1$ has little work per rank to share, as explained in Section~\ref{sec:why}, so that the one-dimensional problems show the cost of the predictor rather than that of the factorization. The runs of Problem~6.1 (Table~\ref{tab:p61}) print a stage profile: for $p=1$ the APSOR stage takes $0.02$ of $0.03$~s at $N=1024$, $0.13$ of $0.15$~s at $N=2048$, $0.48$ of $0.51$~s at $N=4096$ and $3.67$ of $3.71$~s ($99\%$) at $N=8192$, while the SPIKE decomposition rounds to $0.00$~s in every one-dimensional run. In one dimension the cost is therefore the cost of the sequential APSOR sweeps, not of the factorization, which is the opposite of the two-dimensional profile of Table~\ref{tab:profile}. The APSOR stage runs on rank~0 only, so that adding processes cannot shorten it.
	
	\textbf{Growth with $N$.} With $p=1$ the time grows by the factors $1.4$, $1.8$, $1.3$, $2.0$, $2.6$, $2.8$, $5.1$, $3.5$ and $7.3$ per doubling of $N$ from $16$ to $8192$. A tridiagonal solve alone would give a factor $2$. The extra growth comes from the number of APSOR sweeps, each costing $O(N)$: the runs of Problem~6.1 needed $207$, $372$, $1293$, $2560$ and $10\,231$ sweeps for $N=512,\dots,8192$ (Table~\ref{tab:p61}), and the time per sweep, $0.10$, $0.19$ and $0.36$~ms at $N=2048$, $4096$ and $8192$, doubles with $N$ as an $O(N)$ sweep should.
	
	\textbf{Balance.} The ranks finish close together: for $N\ge2048$ the workers agree to within $0.3\%$ and rank~0, which carries APSOR and the reduced solve, is behind them by $0.3\%$ to $5\%$ (for example $5.896$~s against $5.871$~s at $N=8192$, $p=16$).

	\section{Numerical Problems}\label{sec:problems}
	
	We solve eight problems. Four are classical obstacle problems ($g=0$, $f=0$) from Tran et al.\ \cite{TSFO14} and four are soft-layer problems from Cheng et al.\ \cite{CRWX21}. They are stated below and summarized, with the data, the discretization and the reference used, in Table~\ref{tab:problems}. The same solver (Algorithm~\ref{alg:main}) is used for all of them; only the obstacle, the load, $g$ and the domain change. The grid size $N$ is a run-time argument of every code.
	
	\begin{table}[!htbp]
		\centering\small
		\caption{Summary of the eight numerical problems.}\label{tab:problems}
		\begin{tabular}{p{1.9cm}p{2.1cm}p{4.4cm}p{1.9cm}p{0.9cm}p{1.6cm}}
			\toprule
			Problem & Domain & Obstacle & Load $f$ & $g$ & Reference \\
			\midrule
			6.1 & $(0,1)$ & $\varphi_1$ of \eqref{eq:obs12} & $0$ & $0$ & exact \\
			6.2 & $(0,1)$ & $\varphi_2$ of \eqref{eq:obs12} & $0$ & $0$ & exact \\
			6.3 & $[-2,2]^2$ & $\sqrt{1-x^2-y^2}$ (radial) & $0$ & $0$ & analytic \\
			6.4 & $(0,1)^2$ & non-smooth $\varphi_3$ & $0$ & $0$ & none \\
			Example 1 & $(0,1)^2$ & $-\operatorname{dist}((x,y),\partial\Omega)$ & $-8$ & $0.1$ & none \\
			Example 2 & $(0,1)^2$ & $-\operatorname{dist}((x,y),\partial\Omega)$ & piecewise & $0.1$ & none \\
			Example 3 & $(0,1)^2$ & $-\operatorname{dist}((x,y),\partial\Omega)$ & $11(x+y-1)$ & $0.001$ & none \\
			Example 4 & $(0,1)^2$ & $10x(1-x)y(1-y)-\tfrac12$ & $-20$ & $0.1$ & none \\
			\bottomrule
		\end{tabular}
	\end{table}
	
	\subsection*{Classical obstacle problems ($g=0$)}
	
	\paragraph{Problem 6.1: one-dimensional obstacle $\varphi_1$ \cite[Eq.~(6.1)]{TSFO14}.} Find $u=\arg\min_{v\ge\varphi_1}\tfrac12\int_0^1|v'|^2dx$ with $v(0)=v(1)=0$, where $\varphi_1$ is given in \eqref{eq:obs12} and $\varphi_1(x)=\varphi_1(1-x)$ on $[0.5,1]$. The exact solution is $u_1$ of \eqref{eq:exact12}; the contact set is $[\tfrac1{2\sqrt2},1-\tfrac1{2\sqrt2}]\approx[0.354,0.646]$.
	
	\paragraph{Problem 6.2: one-dimensional obstacle $\varphi_2$ \cite[Eq.~(6.2)]{TSFO14}.} The same problem with the trigonometric obstacle $\varphi_2$ of \eqref{eq:obs12}. The exact solution is $u_2$ of \eqref{eq:exact12}, extended by $u(x)=u(1-x)$; the contact set is $[0,0.25]\cup[0.75,1]$.
	
	\paragraph{Problem 6.3: two-dimensional radial obstacle.} On $[-2,2]^2$ the obstacle is $\varphi=\sqrt{1-x^2-y^2}$ for $x^2+y^2\le1$ and $-1$ elsewhere. With $r=\sqrt{x^2+y^2}$ the analytic solution is $u=\sqrt{1-r^2}$ for $r\le r^*$ and $u=-(r^*)^2\log(r/2)/\sqrt{1-(r^*)^2}$ for $r\ge r^*$, where $r^*$ solves $(r^*)^2\big(1-\log(r^*/2)\big)=1$. The code finds $r^*=0.69797$ by bisection, and the boundary condition of the square is the analytic solution itself.
	
	\paragraph{Problem 6.4: two-dimensional non-smooth obstacle \cite[Eq.~(6.4)]{TSFO14}.} On $(0,1)^2$ with $u=0$ on the boundary,
	\begin{equation}\label{eq:phi3}
		\varphi_3(x,y)=\begin{cases}5.0, & |x-0.6|+|y-0.6|<0.04,\\ 4.5, & (x-0.6)^2+(y-0.25)^2<0.001,\\ 4.5, & y=0.57\ \text{and}\ 0.075<x<0.13,\\ 0, & \text{otherwise.}\end{cases}
	\end{equation}
	On the grid, the segment is the row of nodes nearest to $y=0.57$. No exact solution is known, so the tests are the KKT residual and the shape of the solution.
	
	\subsection*{Soft-layer obstacle problems ($g>0$)}
	All four examples use the unit square, $u=0$ on the boundary, stiffness $c_\nu=1$ and the admissible set $u\ge\psi-g$ of Section~\ref{sec:model}. Examples 1 to 3 use the obstacle $\psi=-\operatorname{dist}((x,y),\partial\Omega)$, the distance to the boundary of the square.
	
	\paragraph{Example 1.} Modified from Brezis and Sibony \cite{BS71}: $f=-8$ and $g=0.1$.
	
	\paragraph{Example 2.} A unilateral problem from the literature with a soft layer $g=0.1$ added. The load is
	\[
	f(x,y)=\begin{cases}300, & (x,y)\in S=\{|x-y|\le0.1,\ x\le0.3\},\\ -70\,e^{y}h(x), & x\le1-y,\ (x,y)\notin S,\\ 15\,e^{y}h(x), & x>1-y,\ (x,y)\notin S,\end{cases}
	\]
	where
	\[
	h(x)=\begin{cases}6x, & 0<x\le\tfrac16,\\ 2(1-3x), & \tfrac16<x\le\tfrac13,\\ 6(x-\tfrac13), & \tfrac13<x\le\tfrac12,\\ 2\big(1-3(x-\tfrac13)\big), & \tfrac12<x\le\tfrac23,\\ 6(x-\tfrac23), & \tfrac23<x\le\tfrac56,\\ 2\big(1-3(x-\tfrac23)\big), & \tfrac56<x\le1.\end{cases}
	\]
	
	\paragraph{Example 3.} $f=11(x+y-1)$ and a very thin soft layer, $g=0.001$.
	
	\paragraph{Example 4.} The polynomial obstacle $\psi=10x(1-x)y(1-y)-\tfrac12$, the load $f=-20$ and $g=0.1$.
	
	For these examples the computed regions are classified node by node as rigid contact ($u=\psi-g$), soft-layer contact ($\psi-g<u\le\psi$) and no contact ($u>\psi$).

	\subsection{Timing experiments}
	Separately from the eight problems, the MPI programs were timed for the two-dimensional soft-layer problem (Example 1, $g=0.1$) on grids $N\times N$ with $N=16,32,64,128,256$ and for the one-dimensional obstacle problem (Problem~6.1, matrix orders $N=16,\dots,8192$), on $p=1,2,4,8,16$ processes, with $p\le(N-1)/2$ in two dimensions. The CPU time of rank~0 is reported; rank~0 carries the APSOR stage and the reduced solve, so that it is the slowest rank. The results are in Section~\ref{sec:perf}.
	
	\section{Numerical Results}\label{sec:results}
	
	\subsection{Experimental setup}
	The solver of Sections~\ref{sec:method}--\ref{sec:spike} is implemented in Fortran 90 with MPI in double precision, one code per problem. The tolerances are those of Algorithm~\ref{alg:main} and Section~\ref{sec:how}: an APSOR increment of $10^{-3}$ before the first SPIKE solve, a Picard change of $10^{-11}$ with at most 40 passes, and the KKT tolerances $10^{-10}$, $10^{-8}$ and $10^{-6}$. Accuracy is measured against the exact solutions (Problems 6.1, 6.2) and the analytic solution (Problem 6.3); Problem 6.4 and the four soft-layer examples have no closed form and are judged by the KKT test and by the contact structure. Every run reported in this section used 16 MPI ranks on a laptop (Lenovo IdeaPad Slim 5 16IRL8). The grids are $N=1024$ (1023 unknowns) for Problem~6.1, $N=256$ (255 unknowns) for Problem~6.2, $N=64$ on $[-2,2]^2$ (3969 unknowns, $h=1/16$) for Problem~6.3, $N=256$ (65\,025 unknowns) for Problem~6.4, and $N=64$ (3969 unknowns) for the four soft-layer examples. Each code prints, at the end of the run, the number of APSOR sweeps, the last APSOR increment, the number of SPIKE factorizations (one per active set, i.e.\ one per outer iteration; the Picard passes inside a factorization are not counted), the KKT quantities, the range of the APSOR relaxation parameter $\omega$, the number of rigid-contact nodes, the error against the reference where there is one, and the time of rank~0. Table~\ref{tab:struct} summarizes what the plots show, and Table~\ref{tab:runs} gives the numbers printed by the codes.
	
	\begin{table}[!htbp]
		\centering\small
		\caption{Summary of the computed solutions.}\label{tab:struct}
		\begin{tabular}{p{2.2cm}p{3.6cm}p{7.2cm}}
			\toprule
			Problem & Reference & Observed structure \\
			\midrule
			6.1 & exact solution & contact on $[0.3535,0.6465]$; exact set $[0.35355,0.64645]$ \\
			6.2 & exact solution & contact on $[0,0.25]\cup[0.75,1]$, plateau $10$ in between \\
			6.3 & analytic solution & maximum difference $5.99\times10^{-4}$ \\
			6.4 & none & three peaks above the three obstacle pieces \\
			Example 1 & none & soft-layer contact over most of the square, four rigid lobes \\
			Example 2 & none & one peak; two contact regions (upper left, lower right) \\
			Example 3 & none & contact only near the left and bottom edges \\
			Example 4 & none & rigid contact in the middle inside a soft-layer band \\
			\bottomrule
		\end{tabular}
	\end{table}
	
	\begin{table}[!htbp]
		\centering\small
		\caption{Output of the eight codes, 16 MPI ranks. Sweeps: APSOR sweeps; fact.: SPIKE factorizations (outer iterations); $\omega$: range of the relaxation parameter in APSOR; $|G|_{\rm free}$: largest $|G|$ over free nodes; $G_{\rm act}$: smallest $G$ over active nodes; contact: rigid-contact nodes; $t$: time of rank~0 in seconds.}\label{tab:runs}
		\resizebox{\textwidth}{!}{%
			\begin{tabular}{lrrrcccrr}
				\toprule
				Problem & $N$ & sweeps & fact. & $\omega$ & $|G|_{\rm free}$ & $G_{\rm act}$ & contact & $t$ (s) \\
				\midrule
				6.1 & 1024 & 372 & 6 & 0.761--1.990 & $5.6\times10^{-9}$ & $92.3$ & 301 & 0.063 \\
				6.2 & 256 & 155 & 3 & 1.000--1.989 & $4.7\times10^{-10}$ & $9.69$ & 128 & 0.0032 \\
				6.3 & 64 & 51 & 2 & 0.839--1.919 & $5.1\times10^{-13}$ & $0.59$ & 421 & 0.105 \\
				6.4 & 256 & 189 & 1 & 0.761--1.990 & $4.5\times10^{-9}$ & $0.00$ & 438 & 7.39 \\
				Example 1 & 64 & 37 & 3 & 1.006--1.886 & $6.2\times10^{-12}$ & $0.0835$ & 568 & 0.104 \\
				Example 2 & 64 & 56 & 2 & 1.070--1.933 & $2.1\times10^{-11}$ & $0.194$ & 413 & 0.072 \\
				Example 3 & 64 & 7 & 4 & 0.691--1.167 & $2.3\times10^{-12}$ & $0.453$ & 68 & 0.109 \\
				Example 4 & 64 & 31 & 2 & 1.167--1.869 & $5.7\times10^{-12}$ & $2.56$ & 1305 & 0.097 \\
				\bottomrule
		\end{tabular}}
	\end{table}
	
	\begin{table}[!htbp]
		\centering\small
		\caption{Breakdown of the run time of rank~0 (seconds, printed to two decimals; the total also contains the setup).}\label{tab:runtime}
		\begin{tabular}{lrrrr}
			\toprule
			Problem & APSOR & SPIKE decomposition & Picard / reduced & total \\
			\midrule
			6.1 & 0.05 & 0.00 & 0.00 & 0.06 \\
			6.2 & 0.00 & 0.00 & 0.00 & 0.00 \\
			6.3 & 0.02 & 0.05 & 0.00 & 0.10 \\
			6.4 & 0.88 & 6.35 & 0.07 & 7.39 \\
			Example 1 & 0.01 & 0.07 & 0.01 & 0.10 \\
			Example 2 & 0.01 & 0.04 & 0.00 & 0.07 \\
			Example 3 & 0.00 & 0.09 & 0.01 & 0.11 \\
			Example 4 & 0.02 & 0.06 & 0.01 & 0.10 \\
			\bottomrule
		\end{tabular}
	\end{table}
	
	\subsection{Problems 6.1 and 6.2: one-dimensional obstacles}
\begin{figure}[t]
	\centering
	\includegraphics[width=0.9\linewidth]{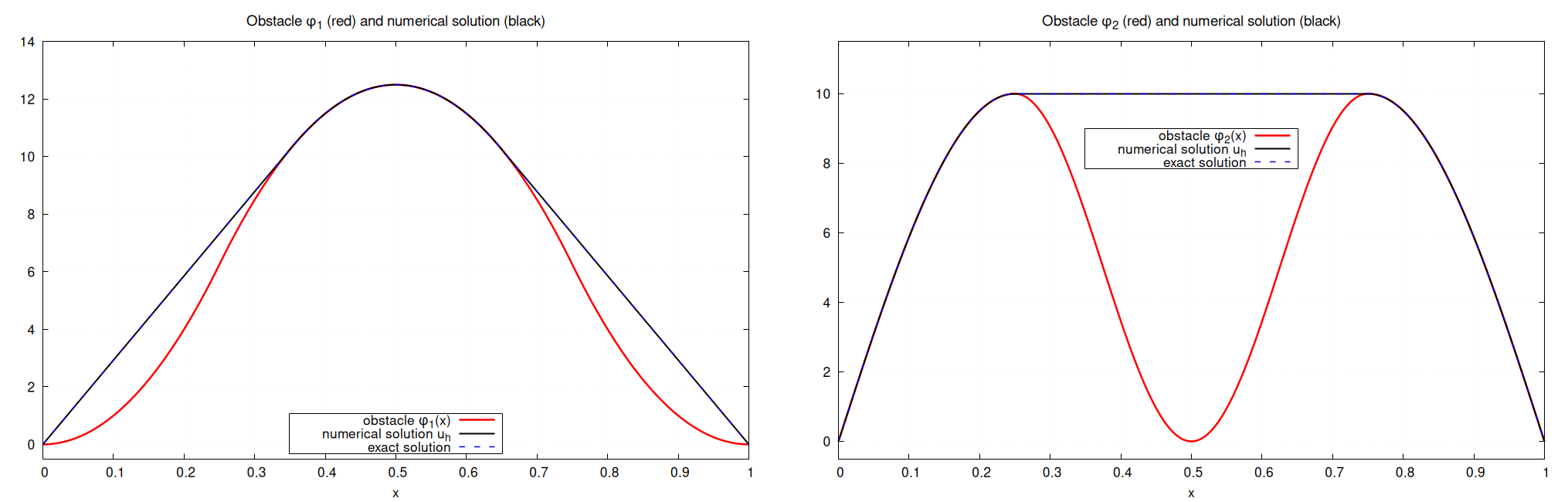}
	\caption{One-dimensional obstacle problems: numerical and exact solutions coincide.}
	\label{fig:figure05problems6162solutions}
\end{figure}

	Figure~\ref{fig:figure05problems6162solutions} shows the computed solutions. In both cases the numerical solution coincides with the exact one of \eqref{eq:exact12}. For $\varphi_1$ it follows the straight line of slope $100-50\sqrt2\approx29.3$ away from the middle and lies on the obstacle over the middle interval, about $[0.35,0.65]$, in agreement with the exact contact set $[0.354,0.646]$. For $\varphi_2$ it follows the obstacle on the two outer arcs $[0,0.25]$ and $[0.75,1]$ and stays at the plateau $10$ in between, where the obstacle dips to $0$ at $x=0.5$.
	
	\paragraph{What the numbers show.} For Problem~6.1 ($N=1024$) the maximum relative error against the exact solution is $1.14\times10^{-8}$, the computed contact set is $[0.35352,0.64648]$ against the exact $[0.35355,0.64645]$, and $301$ nodes rest on the obstacle, which is $0.293\cdot1024$ as the exact contact length $1-1/\sqrt2$ predicts. APSOR needed $372$ sweeps and the active set was frozen and solved with SPIKE six times. For Problem~6.2 ($N=256$) the error is $7.2\times10^{-14}$, the contact set is $[0,0.25]\cup[0.75,1]$ as exactly, and the maximum $10.00000$ is the plateau. Here $128$ nodes are in contact, namely $64$ on each outer arc. The error is at rounding level because the kink of the exact solution at $x=0.25$ is a grid node ($64/256$) and the plateau is constant, so that the discrete solution coincides with the exact one at the nodes; for $\varphi_1$ the kink $x=1/(2\sqrt2)$ is not a node and the error is $1.1\times10^{-8}$ on this fine grid (Section~\ref{sec:pindep} gives the errors on coarser grids). Only three SPIKE factorizations were needed for Problem~6.2 and six for Problem~6.1. In both runs the KKT test passed at its first attempt, with a margin: $|G|\le5.6\times10^{-9}$ on free nodes against the tolerance $10^{-6}$, and $G\ge9.7$ on the active nodes, so that strict complementarity holds with a wide gap.
	
	\subsection{Problem 6.3: radial obstacle}
\begin{figure}[H]
	\centering
	\includegraphics[width=0.9\linewidth]{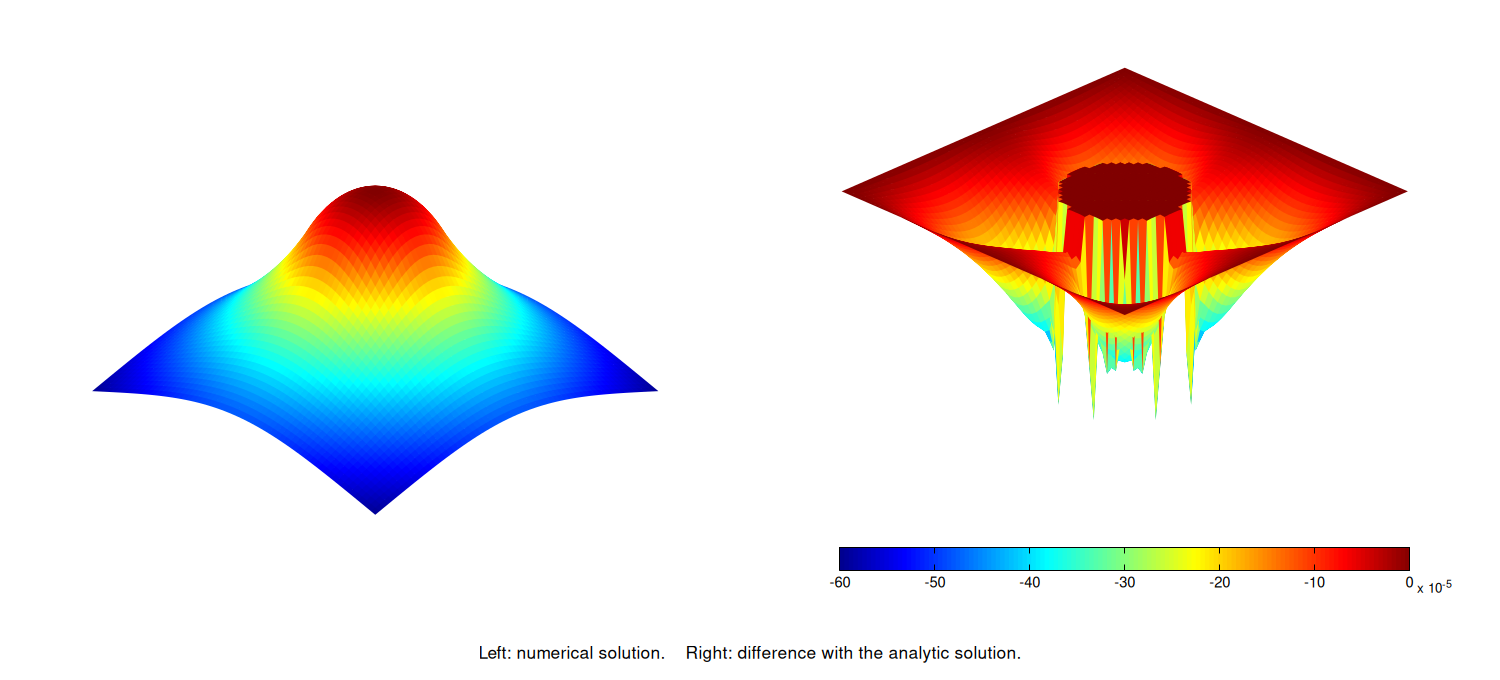}
	\caption{Problem 6.3: numerical solution (left) and its difference from the analytic solution (right), on a colour scale from $-60\times10^{-5}$ to $0$.}
	\label{fig:figure06problem63solutionanderror}
\end{figure}

	Figure~\ref{fig:figure06problem63solutionanderror} shows the numerical solution on the left and, on the right, its difference from the analytic solution on a colour scale from $-60\times10^{-5}$ to $0$. The difference is therefore at most about $6\times10^{-4}$ on the grid of the figure.
	
	\paragraph{What the numbers show.} On the $64\times64$ grid over $[-2,2]^2$ ($h=1/16$) the code reports a maximum error of $5.99\times10^{-4}$ against the analytic solution; it is the absolute error, and the relative error is the same because the maximum of the solution is $1$. The contact radius found on the grid is $0.71261$ against the exact $r^*=0.69797$; the difference, $0.0146$, is a quarter of a mesh width, as it should be for a radius that is read off the grid nodes. The solution has minimum $-0.2142$ and maximum $1.0000$; $421$ nodes are in rigid contact. The solve needed $51$ APSOR sweeps and two SPIKE factorizations. The bandwidth is now $b=63$, and the SPIKE decomposition is $0.05$ of the $0.10$ seconds.
	
	\subsection{Problem 6.4: non-smooth obstacle}
\begin{figure}[H]
	\centering
	\includegraphics[width=0.9\linewidth]{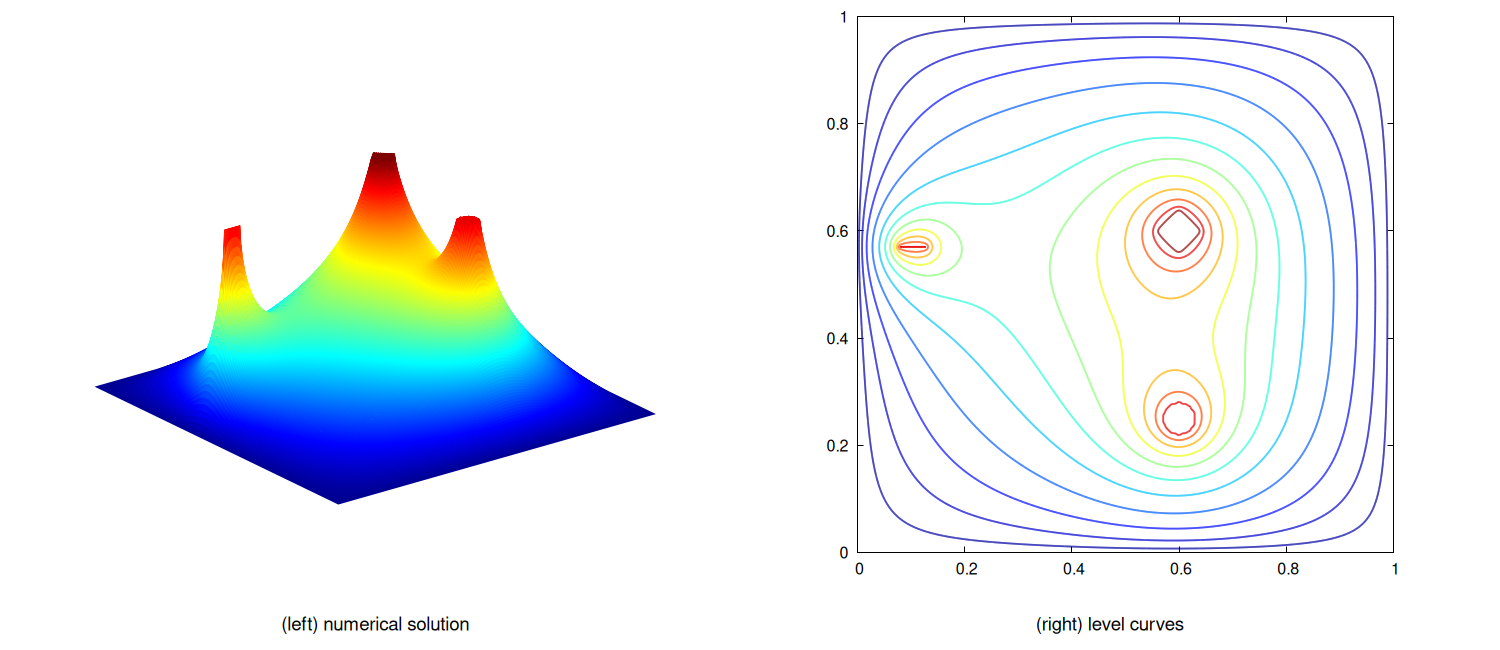}
	\caption{Problem 6.4: numerical solution (left) and level curves (right). The three peaks lie above the diamond, the disk and the segment of the obstacle.}
	\label{fig:figure07problem64solutionandlevelcurves}
\end{figure}

	In Figure~\ref{fig:figure07problem64solutionandlevelcurves} the solution has three peaks above the three obstacle pieces, the diamond and the disk near $x=0.6$ and the segment near $x=0.1$, and it is flat and equal to $0$ elsewhere; the level curves close around each peak. For the $256\times256$ grid the run of the code of this paper (16 ranks) ended after $189$ APSOR sweeps and \emph{one} SPIKE factorization, with the largest and smallest values of the solution $5.00000$ and $0.0003$, $438$ nodes in rigid contact, and $|G|_{\rm free}\le4.5\times10^{-9}$. The peak value $5$ is the height of the diamond piece of the obstacle. There is no closed-form solution, so the correctness of this run rests on the KKT test, which was passed. The SPIKE decomposition is $6.35$ of the $7.39$ seconds ($86\%$) and APSOR $0.88$ seconds ($12\%$). The smallest $G$ on the active nodes is $0$ to the printed precision: some active node carries no contact force, so that the problem is not strictly complementary, and Theorem~\ref{thm:outer} does not apply to this run. The test that was passed is the KKT test itself, not the assumption of the theorem.
	
	\subsection{Examples 1 to 4: soft-layer obstacles}
	In all four examples the computed regions are consistent with the three-region structure of Section~\ref{sec:model}: in the contact plots black marks the nodes with no contact ($u>\psi$), grey the nodes in soft-layer contact ($\psi-g<u<\psi$) and white the nodes in rigid contact ($u=\psi-g$). Table~\ref{tab:runs} shows that the KKT test passed in each example with $G\ge0.08$ on the active nodes, so that strict complementarity holds in all four and Theorem~\ref{thm:outer} applies. The printed minimum and maximum of the solutions are $-0.490$ and $-0.0046$ (Example~1), $-0.513$ and $2.031$ (Example~2), $-0.132$ and $0.134$ (Example~3), and $-0.284$ and $0.025$ (Example~4).
	
\begin{figure}[H]
	\centering
	\includegraphics[width=0.9\linewidth]{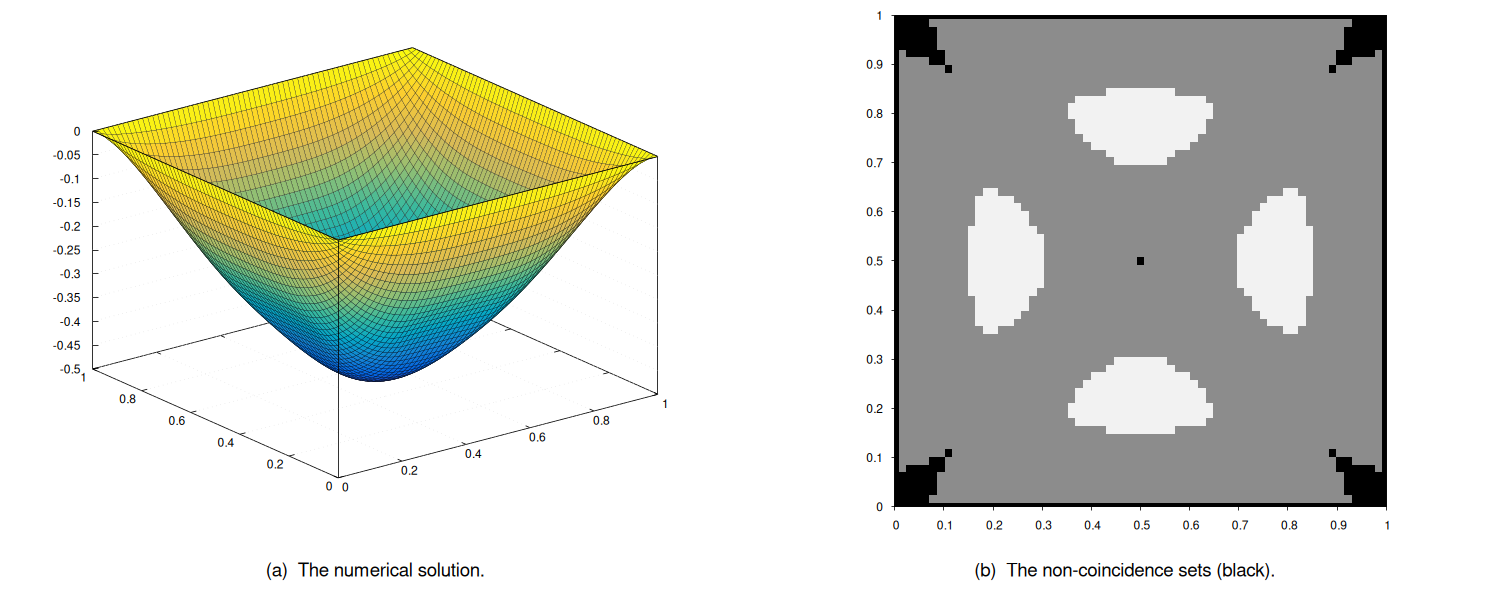}
	\caption{Example 1: numerical solution (left) and contact regions (right). In the right panel black marks the nodes with no contact, grey the nodes in soft-layer contact and white the nodes in rigid contact.}
	\label{fig:figure08example1solutionandcontactregions}
\end{figure}

	\paragraph{Example 1} (Figure~\ref{fig:figure08example1solutionandcontactregions}). The solution is negative everywhere (maximum $-0.0046$) and reaches $-0.490$ at the centre; $37$ APSOR sweeps and three factorizations were needed, and $568$ nodes are in rigid contact. Most of the square is in soft-layer contact, with four rigid-contact lobes near the middle of the sides and no contact in the four corners and at a single node at the centre.
	
\begin{figure}[H]
	\centering
	\includegraphics[width=0.9\linewidth]{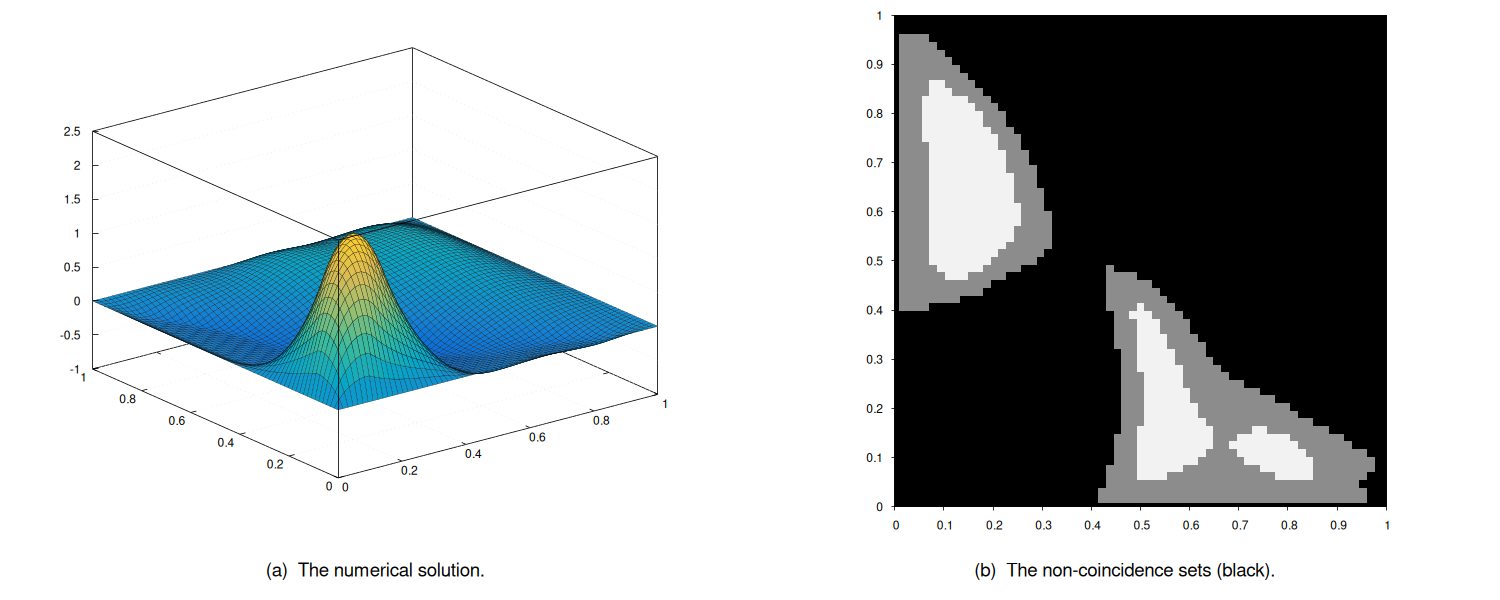}
	\caption{Example 2: numerical solution (left) and contact regions (right). In the right panel black marks the nodes with no contact, grey the nodes in soft-layer contact and white the nodes in rigid contact.}
	\label{fig:figure09example2solutionandcontactregions}
\end{figure}

	\paragraph{Example 2} (Figure~\ref{fig:figure09example2solutionandcontactregions}). The solution (maximum $2.031$) has a single peak and most of the square is not in contact. Soft-layer and rigid contact occur in two regions, one at the upper left and one at the lower right; $413$ nodes are in rigid contact, with $56$ sweeps and two factorizations.
	
\begin{figure}[H]
	\centering
	\includegraphics[width=0.9\linewidth]{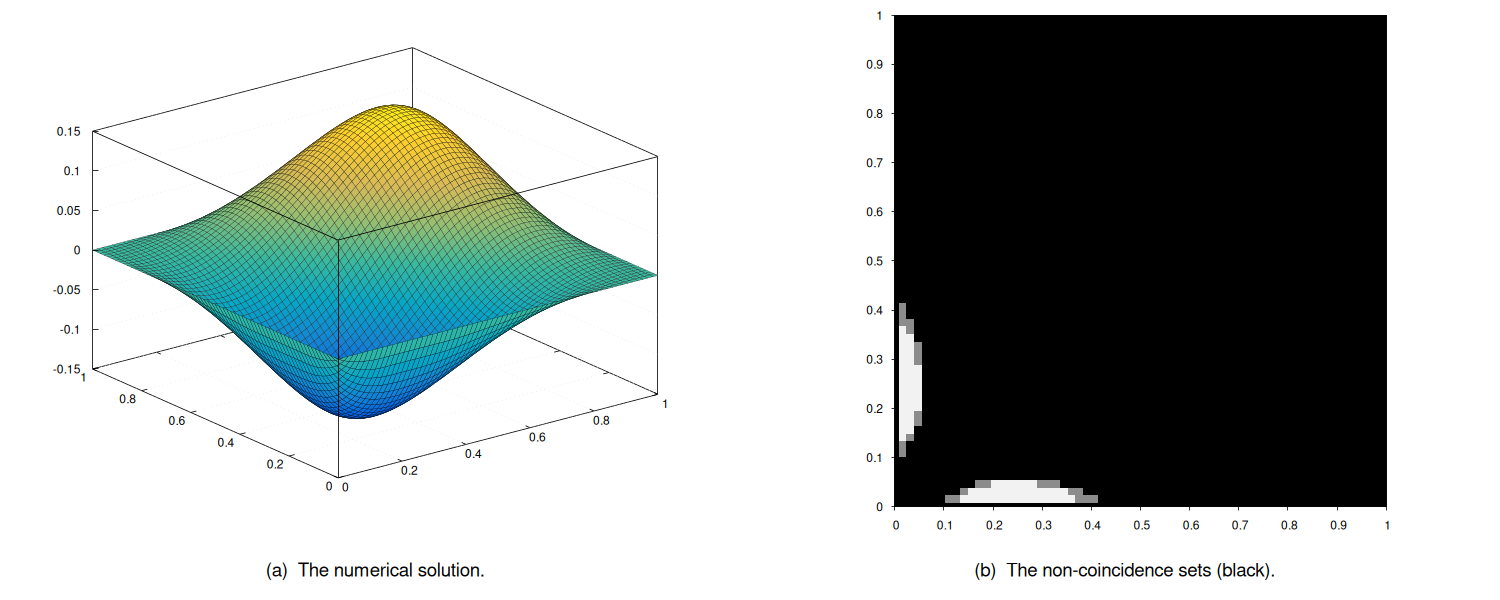}
	\caption{Example 3: numerical solution (left) and contact regions (right). In the right panel black marks the nodes with no contact, grey the nodes in soft-layer contact and white the nodes in rigid contact.}
	\label{fig:figure10example3solutionandcontactregions}
\end{figure}

	\paragraph{Example 3} (Figure~\ref{fig:figure10example3solutionandcontactregions}). The solution ranges between about $-0.15$ and $0.15$. With the very thin layer $g=0.001$ contact occurs only in two small regions near the left edge and the bottom edge, and there is little soft-layer contact. Only $68$ nodes are in rigid contact, and the run needed just $7$ APSOR sweeps (with $\omega$ never above $1.17$) but four factorizations, the most of the eight problems, so that the decomposition is $0.09$ of the $0.11$ seconds.
	
\begin{figure}[H]
	\centering
	\includegraphics[width=0.9\linewidth]{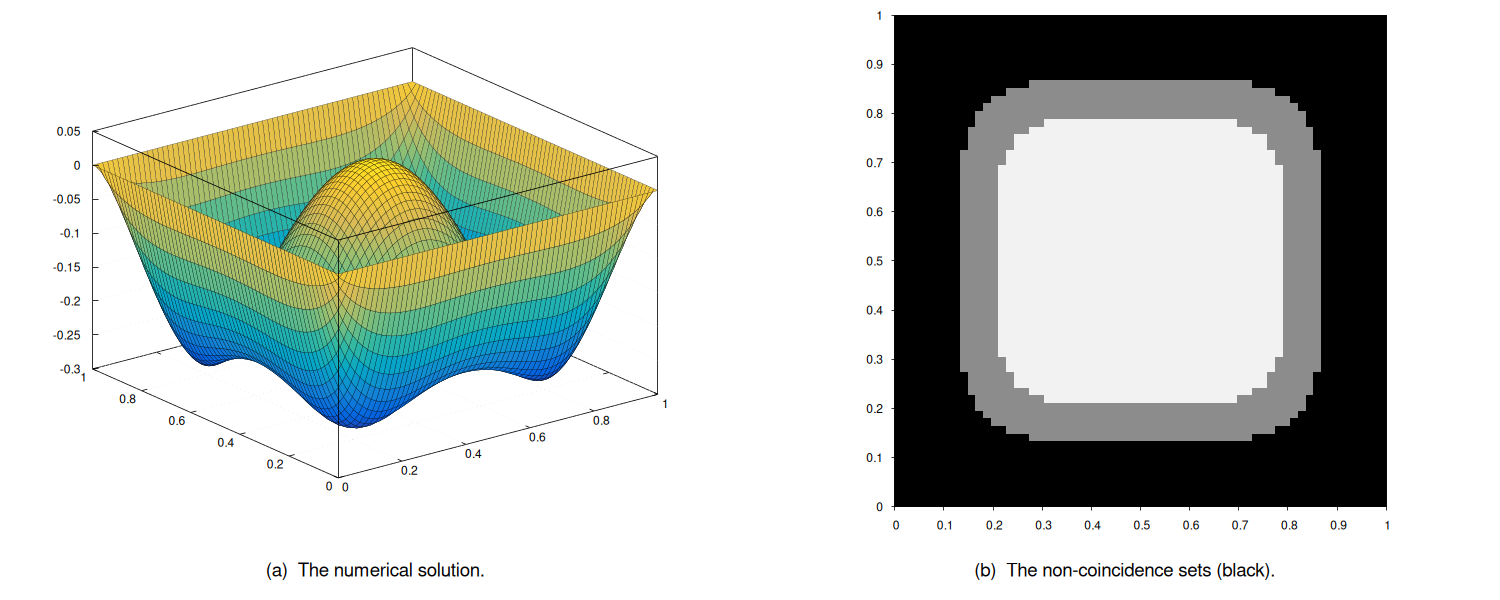}
	\caption{Example 4: numerical solution (left) and contact regions (right). In the right panel black marks the nodes with no contact, grey the nodes in soft-layer contact and white the nodes in rigid contact.}
	\label{fig:figure11example4solutionandcontactregions}
\end{figure}

	\paragraph{Example 4} (Figure~\ref{fig:figure11example4solutionandcontactregions}). The solution ranges between about $-0.3$ and $0.05$. A large rigid-contact region in the middle is surrounded by a band of soft-layer contact, and the rest of the square is not in contact. The rigid-contact region is the largest of the four, $1305$ of $3969$ nodes, and the active forces are the largest, $G\ge2.56$; $31$ sweeps and two factorizations were needed.
	
	\subsection{Problem 6.1 on ten grids and 1 to 16 processes}\label{sec:pindep}
	Problem~6.1 was run on the grids $N=16,32,\dots,8192$ with up to $16$ ranks ($p=1,2,4$ for $N=16$; $p=1,2,4,8$ for $N=32$; $p=1,2,4,8,16$ for $N\ge64$). The one-dimensional timing study of Section~\ref{sec:perf1d} was made with these runs. Table~\ref{tab:p61} gives what the code prints.
	
	\begin{table}[!htbp]
		\centering\small
		\caption{Problem 6.1: output of the code on ten grids. For every $N$ the number of sweeps, the number of factorizations, $G_{\rm act}$, the number of contact nodes, the error and the contact set are identical for all the $p$ run. The last column gives the smallest and largest value of $|G|_{\rm free}$ over the $p$ run.}\label{tab:p61}
		\resizebox{\textwidth}{!}{%
			\begin{tabular}{rcrrrrrcc}
				\toprule
				$N$ & $p$ run & sweeps & fact. & $G_{\rm act}$ & contact & error & contact set & $|G|_{\rm free}$ \\
				\midrule
				16 & 1,2,4 & 19 & 1 & 167 & 5 & $3.0664\times10^{-3}$ & $[0.37500,0.62500]$ & $9.09\times10^{-13}$ \\
				32 & 1--8 & 21 & 1 & 36.4 & 11 & $7.6885\times10^{-4}$ & $[0.34375,0.65625]$ & $5.46\times10^{-12}$ \\
				64 & 1--16 & 38 & 2 & 174 & 19 & $2.5934\times10^{-4}$ & $[0.35938,0.64062]$ & $[1.46,\,2.18]\times10^{-11}$ \\
				128 & 1--16 & 55 & 2 & 48.9 & 39 & $3.1709\times10^{-5}$ & $[0.35156,0.64844]$ & $[5.82,\,11.6]\times10^{-11}$ \\
				256 & 1--16 & 100 & 3 & 198 & 75 & $2.9026\times10^{-5}$ & $[0.35547,0.64453]$ & $[2.33,\,3.49]\times10^{-10}$ \\
				512 & 1--16 & 207 & 5 & 96.1 & 151 & $1.1410\times10^{-8}$ & $[0.35352,0.64648]$ & $1.40\times10^{-9}$ \\
				1024 & 1--16 & 372 & 6 & 92.3 & 301 & $1.1410\times10^{-8}$ & $[0.35352,0.64648]$ & $5.59\times10^{-9}$ \\
				2048 & 1--16 & 1293 & 8 & 84.5 & 601 & $1.1410\times10^{-8}$ & $[0.35352,0.64648]$ & $2.24\times10^{-8}$ \\
				4096 & 1--16 & 2560 & 9 & 69.1 & 1201 & $1.1410\times10^{-8}$ & $[0.35352,0.64648]$ & $[8.94,\,11.9]\times10^{-8}$ \\
				8192 & 1--16 & 10231 & 11 & 38.1 & 2401 & $1.1410\times10^{-8}$ & $[0.35352,0.64648]$ & $[3.58,\,4.77]\times10^{-7}$ \\
				\bottomrule
		\end{tabular}}
	\end{table}
	
	\paragraph{Independence of the number of processes.} Proposition~\ref{prop:spike} states that the SPIKE solve gives the same result for every number of processes. In all $47$ runs the printed solution quantities agree for the different $p$ on the same grid: the number of APSOR sweeps, the number of factorizations, the minimum and maximum of $u$, the number of contact nodes, the error and the contact set. Three printed quantities differ slightly. (a) The KKT residual $|G|_{\rm free}$ is identical for all $p$ for $N\le32$ and for $N=512,1024,2048$, and differs for $N=64$ (at $p=16$), $N=128$, $N=256$ (at $p=8$), $N=4096$ (at $p=8$) and $N=8192$ (at $p=4$), always by a factor of at most about $2$ and always far below the tolerance $10^{-6}$. (b) The last APSOR increment differs for $N=2048$, $4096$ and $8192$ (for example $1.3\times10^{-11}$, $1.0\times10^{-10}$ and $2.4\times10^{-10}$ at $N=2048$ for $p=1$, $4$ and $8$), and (c) the smallest relaxation parameter in APSOR differs for $N=2048$, $4096$ and $8192$ (at $N=2048$: $0.839$ for $p=1,2,8$ and $0.686$ for $p=4,16$) and the largest for $N=256$ ($1.988$ at $p=4$, $1.986$ otherwise).
	
	\paragraph{The error against the exact solution.} The error is not a smooth function of $N$ (Table~\ref{tab:p61}): it is $3.07\times10^{-3}$, $7.69\times10^{-4}$, $2.59\times10^{-4}$, $3.17\times10^{-5}$, $2.90\times10^{-5}$ for $N=16,\dots,256$ and then exactly $1.1410\times10^{-8}$ for every $N$ from $512$ to $8192$. The reason is the kink of the exact solution at $x^*=1/(2\sqrt2)=0.353553$, where $u''$ jumps from $0$ on the line to $200$ on the obstacle $100x^2$. The node nearest to $x^*$ is at distance $d$, and the grids $N=512,\dots,8192$ all contain the same node $181/512=0.353516$ at $d=3.78\times10^{-5}$. A Taylor estimate gives an error $\approx100\,d^2$ at that node; the code prints the error relative to $\max u=12.5$, so the prediction is $100\,d^2/12.5$. Table~\ref{tab:err61} compares it with the output. This is our estimate and not a theorem; it agrees to four or five digits for $N=32$, $128$ and $512$ to $8192$, to within $5\%$ for $N=64$ and $256$ and within $20\%$ for $N=16$. The error is therefore bounded by a multiple of $h^2$ but is controlled by the position of the kink relative to the grid, and we do not claim a convergence rate from these runs.
	
	\begin{table}[!htbp]
		\centering\small
		\caption{Problem 6.1: distance $d$ from the kink $x^*=1/(2\sqrt2)$ to the nearest node, the estimate $100\,d^2/12.5$ and the error printed by the code.}\label{tab:err61}
		\begin{tabular}{rccc}
			\toprule
			$N$ & $d$ & $100\,d^2/12.5$ & error (output) \\
			\midrule
			16 & $2.145\times10^{-2}$ & $3.68\times10^{-3}$ & $3.0664\times10^{-3}$ \\
			32 & $9.80\times10^{-3}$ & $7.689\times10^{-4}$ & $7.6885\times10^{-4}$ \\
			64 & $5.82\times10^{-3}$ & $2.71\times10^{-4}$ & $2.5934\times10^{-4}$ \\
			128 & $1.99\times10^{-3}$ & $3.171\times10^{-5}$ & $3.1709\times10^{-5}$ \\
			256 & $1.92\times10^{-3}$ & $2.94\times10^{-5}$ & $2.9026\times10^{-5}$ \\
			512 to 8192 & $3.78\times10^{-5}$ & $1.141\times10^{-8}$ & $1.1410\times10^{-8}$ \\
			\bottomrule
		\end{tabular}
	\end{table}
	
	\paragraph{The KKT residual grows with $N$.} From $N=512$ the largest $|G|_{\rm free}$ is $1.40\times10^{-9}$, $5.59\times10^{-9}$, $2.24\times10^{-8}$, $8.94\times10^{-8}$ and $3.58\times10^{-7}$ at $N=512,\dots,8192$, a factor $4$ per doubling. This is the scaling of the stiffness entries $h^{-2}$ and suggests that the residual is rounding error amplified by $h^{-2}$, and it stays below the tolerance $10^{-6}$ on all grids tested ($N\le8192$). The count of contact nodes follows the exact contact length: $2401$ nodes at $N=8192$ against $0.2929\cdot8192=2399.6$.
	
	\subsection{Example 1 on five grids and 1 to 16 processes}\label{sec:ex1grids}
	Example 1 ($g=0.1$) was run on $N=16,32,64,128,256$ with $p=1,2,4$ for $N=16$, $p=1,2,4,8$ for $N=32$ and $p=1,2,4,8,16$ for $N\ge64$. The two-dimensional timing study of Section~\ref{sec:perf} was made with these runs. Table~\ref{tab:ex1} gives what the code prints.
	
	\begin{table}[!htbp]
		\centering\small
		\caption{Example 1: output of the code on five grids. For every $N$ the number of sweeps, the number of factorizations, $G_{\rm act}$, the extreme values of $u$, the number of contact nodes and the range of $\omega$ are identical for all the $p$ run. The last column gives the smallest and largest value of $|G|_{\rm free}$ over the $p$ run.}\label{tab:ex1}
		\resizebox{\textwidth}{!}{%
			\begin{tabular}{rrrrrrrrrc}
				\toprule
				$N$ & unknowns & sweeps & fact. & $G_{\rm act}$ & $\min u$ & $\max u$ & contact & $\omega$ & $|G|_{\rm free}$ \\
				\midrule
				16 & 225 & 22 & 2 & 1.12 & $-0.48999$ & $-0.04565$ & 44 & 1.026--1.685 & $[1.29,\,1.32]\times10^{-12}$ \\
				32 & 961 & 33 & 2 & 0.947 & $-0.48977$ & $-0.01486$ & 152 & 1.167--1.888 & $[7.13,\,7.53]\times10^{-12}$ \\
				64 & 3969 & 37 & 3 & 0.0835 & $-0.48963$ & $-0.00458$ & 568 & 1.006--1.886 & $[6.20,\,8.12]\times10^{-12}$ \\
				128 & 16129 & 58 & 4 & 0.123 & $-0.48964$ & $-0.00136$ & 2028 & 1.006--1.961 & $[4.57,\,5.87]\times10^{-11}$ \\
				256 & 65025 & 89 & 5 & 0.333 & $-0.48963$ & $-0.00039$ & 7980 & 1.006--1.985 & $[2.49,\,4.95]\times10^{-10}$ \\
				\bottomrule
		\end{tabular}}
	\end{table}
	
	\paragraph{Independence of the number of processes.} In all $22$ runs the printed solution quantities agree for the different $p$ on the same grid: sweeps, factorizations, $G_{\rm act}$, the extreme values of $u$, the number of contact nodes and the last APSOR increment. The only printed quantity that differs is the KKT residual $|G|_{\rm free}$ (last column of Table~\ref{tab:ex1}, always by a factor of at most $2$ and far below $10^{-6}$), which we read as rounding, as for Problem~6.1. There is no exact solution, so the runs give no error. As for Problem~6.1, the comparison is made on printed quantities and not on the solution vectors.
	
	\paragraph{The solution on finer grids.} The minimum of $u$ stays near $-0.4896$ ($-0.48999$, $-0.48977$, $-0.48963$, $-0.48964$, $-0.48963$), the maximum approaches $0$ from below ($-0.0457$, $-0.0149$, $-0.0046$, $-0.0014$, $-0.0004$, decreasing by factors $3.1$, $3.2$, $3.4$, $3.5$ per doubling), and the share of nodes in rigid contact falls slowly ($19.6\%$, $15.8\%$, $14.3\%$, $12.6\%$, $12.3\%$ of the unknowns). We do not claim a rate, since there is no reference. Both the sweeps ($22$, $33$, $37$, $58$, $89$) and the factorizations ($2$, $2$, $3$, $4$, $5$) grow slowly with $N$, and the active set stays strictly complementary on every grid ($G_{\rm act}\ge0.08$).
	
	\subsection{Problem 6.4 on five grids and 1 to 16 processes}\label{sec:p64grids}
	Problem~6.4 was run on $N=16,32,64,128,256$ with the same process counts as Example~1: $p=1,2,4$ for $N=16$, $p=1,2,4,8$ for $N=32$ and $p=1,2,4,8,16$ for $N\ge64$. Table~\ref{tab:p64out} gives what the code prints. For every $N$ the number of sweeps, the number of factorizations, the number of contact nodes, $G_{\rm act}$, the extreme values of $u$ and the range of $\omega$ are the same for all $p$ run.
	
	\begin{table}[!htbp]
		\centering\small
		\caption{Problem 6.4: output of the code on five grids (identical for all $p$).}\label{tab:p64out}
		\begin{tabular}{rrrrrrrr}
			\toprule
			$N$ & unknowns & sweeps & fact. & contact & $G_{\rm act}$ & $\min u$ & $\omega$ \\
			\midrule
			16 & 225 & 34 & 1 & 2 & $2.0\times10^{3}$ & 0.02364 & 0.761--1.704 \\
			32 & 961 & 49 & 1 & 7 & $2.0\times10^{3}$ & 0.01730 & 0.761--1.877 \\
			64 & 3969 & 78 & 1 & 28 & 0 & 0.00472 & 0.761--1.958 \\
			128 & 16129 & 131 & 1 & 112 & 0 & 0.00120 & 0.761--1.982 \\
			256 & 65025 & 189 & 1 & 438 & 0 & 0.00030 & 0.761--1.990 \\
			\bottomrule
		\end{tabular}
	\end{table}
	
	The maximum of $u$ is $4.5$ for $N=16$ and $5.0$ for all finer grids. The sweeps grow slowly with $N$ ($34$, $49$, $78$, $131$, $189$) and a single SPIKE factorization is enough on every grid.
	
	\begin{table}[!htbp]
		\centering\small
		\caption{Problem 6.4: rank-0 time $T_p$ in seconds.}\label{tab:p64time}
		\begin{tabular}{rrrrrr}
			\toprule
			$N$ & $p=1$ & $p=2$ & $p=4$ & $p=8$ & $p=16$ \\
			\midrule
			16 & 0.0028 & 0.0033 & 0.0024 & -- & -- \\
			32 & 0.0116 & 0.0159 & 0.0160 & 0.0203 & -- \\
			64 & 0.0717 & 0.119 & 0.107 & 0.136 & 0.131 \\
			128 & 0.715 & 1.472 & 1.087 & 1.356 & 1.044 \\
			256 & 8.337 & 23.35 & 17.14 & 21.16 & 13.05 \\
			\bottomrule
		\end{tabular}
	\end{table}
	
	\begin{table}[!htbp]
		\centering\small
		\caption{Problem 6.4: speed-up $S_p=T_1/T_p$ and efficiency $E_p=S_p/p$.}\label{tab:p64sp}
		\resizebox{\textwidth}{!}{%
			\begin{tabular}{r|rr|rr|rr|rr|rr}
				\toprule
				& \multicolumn{2}{c|}{$p=1$} & \multicolumn{2}{c|}{$p=2$} & \multicolumn{2}{c|}{$p=4$} & \multicolumn{2}{c|}{$p=8$} & \multicolumn{2}{c}{$p=16$} \\
				$N$ & $S_p$ & $E_p$ & $S_p$ & $E_p$ & $S_p$ & $E_p$ & $S_p$ & $E_p$ & $S_p$ & $E_p$ \\
				\midrule
				16 & 1.00 & 100\% & 0.86 & 43\% & 1.18 & 29\% & -- & -- & -- & -- \\
				32 & 1.00 & 100\% & 0.73 & 36\% & 0.72 & 18\% & 0.57 & 7\% & -- & -- \\
				64 & 1.00 & 100\% & 0.60 & 30\% & 0.67 & 17\% & 0.53 & 7\% & 0.55 & 3\% \\
				128 & 1.00 & 100\% & 0.49 & 24\% & 0.66 & 16\% & 0.53 & 7\% & 0.68 & 4\% \\
				256 & 1.00 & 100\% & 0.36 & 18\% & 0.49 & 12\% & 0.39 & 5\% & 0.64 & 4\% \\
				\bottomrule
		\end{tabular}}
	\end{table}
	
\begin{figure}[tbph!]
	\centering
	\includegraphics[width=\textwidth]{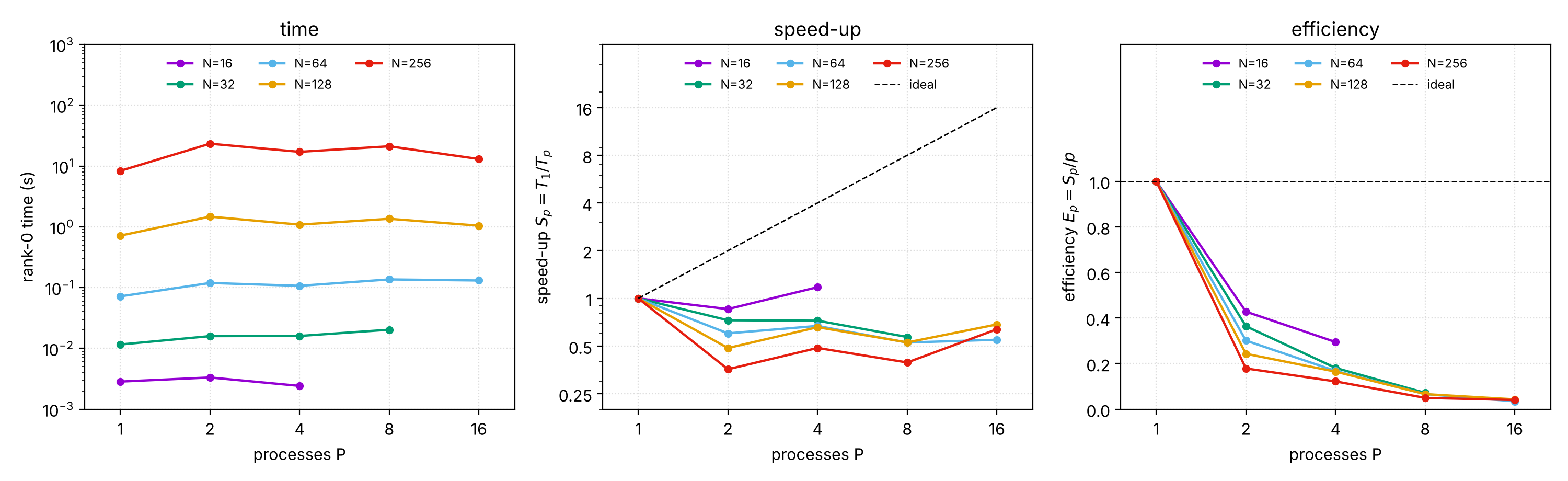}
	\caption{Problem 6.4: rank-0 time, speed-up and efficiency against the number of processes for $N=16,\dots,256$.}
	\label{fig:figure12problem64timespeedupefficiency}
\end{figure}
	
	\paragraph{Timings.} The time grows with $N$ as the size of the banded system does: from $2.8$~ms at $N=16$ to $8.3$~s at $N=256$ on one process. In the $N=256$ run the SPIKE decomposition takes $6.88$~s of $8.34$~s on one process and $11.46$~s of $13.05$~s on $16$ processes, while APSOR needs $0.79$ to $1.33$~s and the Picard stage $0.10$ to $0.52$~s. The largest grid is therefore dominated by the part of the algorithm that is distributed over the ranks, as in Table~\ref{tab:profile}. The time at $p=16$ is $13.05$~s here, against $7.39$~s in the earlier run of Table~\ref{tab:runs}; the two runs have the same iteration counts and the same solution, and the difference is the run-to-run variation of the timing. The measured speed-up (Table~\ref{tab:p64sp}, Figure~\ref{fig:figure12problem64timespeedupefficiency}) is below one in all runs except $N=16$ with $p=4$ ($S_4=1.18$); the efficiency falls from between $18\%$ and $43\%$ at $p=2$ to $3\%$--$4\%$ at $p=16$ for $N\ge64$, and the data of Table~\ref{tab:p64time} are reported as measured.
	
	\subsection{Summary of verification}
	Where an independent reference exists, the computed solution agrees with it: the maximum error is $1.1\times10^{-8}$ for Problem~6.1, $7.2\times10^{-14}$ for Problem~6.2 and $6.0\times10^{-4}$ for Problem~6.3 on a $64\times64$ grid, and the contact sets of Problems 6.1 and 6.2 agree with the exact ones to the width of one cell or better. For Problem~6.4 (also on five grids from $N=16$ to $N=256$ with up to $16$ ranks, Table~\ref{tab:p64out}) and the four soft-layer examples there is no closed form, but the KKT test was passed in all five runs, with $|G|\le5\times10^{-9}$ on free nodes against the tolerance $10^{-6}$ and no negative $G$ on active nodes, and the contact structures are those of the plots. In the eight runs of Table~\ref{tab:runs} the APSOR stage needed from $7$ to $372$ sweeps and the number of SPIKE factorizations was between one and six; Problem~6.1 on the finer grids needed up to $10\,231$ sweeps and $11$ factorizations (Table~\ref{tab:p61}). The timings (Section~\ref{sec:perf}) show that in two dimensions the SPIKE decomposition dominates the cost, while in one dimension the sequential APSOR stage does ($99\%$ at $N=8192$). In two dimensions the decomposition takes $86\%$ of the time for Problem~6.4 and between 0.04 and 0.09 seconds out of the 0.07 to 0.11 seconds of each $64\times64$ run. The independence of the number of processes (Proposition~\ref{prop:spike}) was confirmed for Problem~6.1 on ten grids from $N=16$ to $N=8192$ with up to $16$ ranks (Table~\ref{tab:p61}): the solution quantities are identical, and only the KKT residual, the last APSOR increment and the smallest relaxation parameter differ, at rounding level. The same holds for Example 1 on five grids with up to $16$ ranks (Table~\ref{tab:ex1}).
	
	\section{Summary and conclusions}\label{sec:concl}
	
	We described an APSOR--SPIKE solver for classical and soft-layer obstacle problems. APSOR finds the contact set cheaply, SPIKE solves the resulting banded system exactly and in parallel, a Picard iteration handles the soft-layer compliance with one SPIKE solve per pass, and a KKT test returns the run to APSOR when needed. We proved that APSOR converges, that the SPIKE solve is exact and independent of the number of processes, that the Picard iteration contracts with a factor below $0.052$ on the unit square, and that the loop ends at the unique KKT point under strict complementarity.
	
	On eight test problems the computed solutions agree with the exact or analytic references where these exist, pass the KKT test in every case, and show the expected contact structures, in one and two dimensions and for soft layers from $g=0.001$ to $g=0.1$. The timing study up to $256\times256$ and 16 processes shows where the cost is: in two dimensions the SPIKE decomposition takes $84\%$ to $97\%$ of the run time and APSOR about $1\%$, so the part that is distributed over the ranks is the dominant part of the work. In one dimension (Problem~6.1), matrix orders up to $8192$ show that the cost is that of the sequential APSOR sweeps and not of the factorization. On ten grids the solution did not depend on the number of processes beyond rounding in the KKT residual.
	
	The advantages of the scheme are the following. It needs no parameter tuning, because APSOR chooses its relaxation parameter by the Wolfe conditions. It solves the classical and the soft-layer problem with the same code. The linear system is solved exactly, so that the accuracy reaches $10^{-8}$ to $10^{-14}$ where an exact solution is known, and the result does not depend on the number of processes. The factorization is computed once per active set and reused in every Picard pass, so that only one to six factorizations were needed in the two-dimensional runs. The work and the storage of the factorization per rank decrease with the number of ranks. Every run ends with a KKT certificate, and the convergence of APSOR, of the Picard iteration and of the outer loop is proved.
	
	Further work includes parallelizing the reduced solve or using the recursive and truncated SPIKE variants of \cite{PS06}, and extending the scheme to three dimensions and to other free boundary problems.
	

\end{document}